# ESTIMATORS OF DIFFUSIONS WITH RANDOMLY SPACED DISCRETE OBSERVATIONS: A GENERAL THEORY

By Yacine Aït-Sahalia[1] and Per A. Mykland[2]

*Princeton University and University of Chicago*

We provide a general method to analyze the asymptotic properties of a variety of estimators of continuous time diffusion processes when the data are not only discretely sampled in time but the time separating successive observations may possibly be random. We introduce a new operator, the generalized infinitesimal generator, to obtain Taylor expansions of the asymptotic moments of the estimators. As a special case, our results apply to the situation where the data are discretely sampled at a fixed nonrandom time interval. We include as specific examples estimators based on maximum-likelihood and discrete approximations such as the Euler scheme.

**1. Introduction.** Most theoretical models in finance are spelled out in continuous time [see, e.g., Merton (1992)], whereas the observed data are, by nature, discretely sampled in time. Estimating these models from discrete time observations has become in recent years an active area of research in statistics and econometrics, and a number of estimation procedures have been proposed in the context of parametric models for continuous-time Markov processes, often in the special case of diffusions. Not only are the observations sampled discretely in time, but it is often the case with financial data that the time separating successive observations is itself random, as illustrated, for example, in Figure 1 of Aït-Sahalia and Mykland [(2003), page 484].

This earlier paper focused on the case of inference with the help of likelihood. For data of the type we consider, however, it is common to use a variety of estimating equations, of which likelihood is only one instance; see, for example, Hansen and Scheinkman (1995), Aït-Sahalia (1996, 2002) and

Received June 2002; revised November 2003.
[1]Supported in part by an Alfred P. Sloan Research Fellowship and NSF Grants SBR-99-96023 and SES-01-11140.
[2]Supported in part by NSF Grants DMS-99-71738 and DMS-02-04639.

Bibby, Jacobsen and Sørensen (2004). Our objective in this paper is to carry out a detailed analysis of the asymptotic properties of a large class of such estimators in the context of discretely and randomly sampled data. Unlike Aït-Sahalia and Mykland (2003), it will also permit the diffusion function to depend on both the parameter and the data.

We model this situation as follows. Suppose that we observe the process

$$dX_t = \mu(X_t; \theta) \, dt + \sigma(X_t; \gamma) \, dW_t \tag{1}$$

at discrete times in the interval $[0, T]$, and we wish to estimate the parameters $\theta$ and/or $\gamma$. We call the observation times $\tau_0 = 0, \tau_1, \tau_2, \ldots, \tau_{N_T}$, where $N_T$ is the smallest integer such that $\tau_{N_T+1} > T$. Because the properties of estimators vary widely depending upon whether the drift or the diffusion parameters, or both, are estimated, we consider the three cases of estimating $\beta = (\theta, \gamma)$ jointly, $\beta = \theta$ with $\gamma$ known or $\beta = \gamma$ with $\theta$ known. In regular circumstances, $\hat{\beta}$ converges in probability to some $\bar{\beta}$ and $\sqrt{T}(\hat{\beta} - \bar{\beta})$ converges in law to $N(0, \Omega_\beta)$ as $T$ tends to infinity.

For each estimator, the corresponding $\Omega_\beta$ and, when applicable the bias $\bar{\beta} - \beta_0$, depend on the transition density of the diffusion process, which is generally unknown in closed form. Our solution is to derive Taylor expansions for the asymptotic variance and bias starting with a leading term that corresponds to the limiting case where the sampling is continuous in time. Our main results deliver closed form expressions for the terms of these Taylor expansions. For that purpose, we introduce a new operator, which we call the generalized infinitesimal generator of the diffusion.

Specifically, we write the law of the sampling intervals $\Delta_n = \tau_n - \tau_{n-1}$ as

$$\Delta = \varepsilon \Delta_0, \tag{2}$$

where $\Delta_0$ has a given finite distribution and $\varepsilon$ is deterministic. Our Taylor expansions take the form

$$\Omega_\beta = \Omega_\beta^{(0)} + \varepsilon \Omega_\beta^{(1)} + \varepsilon^2 \Omega_\beta^{(2)} + O(\varepsilon^3), \tag{3}$$

$$\bar{\beta} - \beta_0 = b^{(1)} \varepsilon + b^{(2)} \varepsilon^2 + O(\varepsilon^3). \tag{4}$$

While the limiting term as $\varepsilon$ goes to zero corresponds to continuous sampling, by adding higher-order terms in $\varepsilon$, we progressively correct this leading term for the discreteness of the sampling. The two equations (3) and (4) can then be used to analyze the relative merits of different estimation approaches, by comparing the order in $\varepsilon$ at which various effects manifest themselves, and when they are equal, the relative magnitudes of the corresponding coefficients in the expansion.

Because the coefficients of the expansions depend upon the distribution of the sampling intervals, we can also use these expressions to assess the effect of different sampling patterns on the overall properties of estimators. Moreover, our results apply not only to random sampling, but also to the situation where the sampling interval is time-varying in a deterministic manner (see Section 5.3), or to the case where the sampling interval is simply fixed, in which case we just need to set $\text{Var}[\Delta_0] = 0$ in all our expressions. One particular example is indeed sampling at a deterministic fixed time interval, such as, say, daily or weekly, which is the setup adopted by much of the recent literature on discretely observed diffusions [see, e.g., Hansen and Scheinkman (1995), Aït-Sahalia (1996, 2002) and Bibby, Jacobsen and Sørensen (2004)].

The paper is organized as follows. Section 2 sets up the model and the assumptions used throughout the paper. Section 3 develops a general theory that establishes the asymptotic properties of a large class of estimators of parametric diffusions and their Taylor expansions. Section 4 applies these results to two specific examples of estimating equations: first, the maximum likelihood estimator; second, the Euler approximate discrete scheme based on a Gaussian likelihood. Our conclusions also carry over to the maximum likelihood-type estimators discussed in Aït-Sahalia and Mykland (2003). We discuss extensions of the theory in Section 5. Proofs are contained in Section 6. Section 7 concludes.

## 2. Data structure and inference scheme.

2.1. *The process and the sampling.* We let $\mathcal{S} = (\underline{x}, \bar{x})$ denote the domain of the diffusion $X_t$. In general, $\mathcal{S} = (-\infty, +\infty)$, but in many examples in finance, we are led to consider variables such as nominal interest rates, in which case $\mathcal{S} = (0, +\infty)$. Whenever we are estimating parameters, we will take the parameter space for the $d$-dimensional vector $\beta$ to be an open and bounded set. We will make use of the scale and speed densities of the process, defined as

$$(5) \qquad s(x;\beta) \equiv \exp\left\{-2\int^x (\mu(y;\theta)/\sigma^2(y;\gamma))\,dy\right\},$$

$$(6) \qquad m(x;\beta) \equiv 1/(\sigma^2(x;\gamma)s(x;\beta))$$

and the scale and speed measures $S(x;\beta) \equiv \int^x s(w;\beta)\,dw$ and $M(x;\beta) \equiv \int^x m(w;\beta)\,dw$. The lower bound of integration is an arbitrary point in the interior of $\mathcal{S}$. We also define the same increasing transformation as in Aït-Sahalia (2002):

$$(7) \qquad g(x;\gamma) \equiv \int^x \frac{du}{\sigma(x;\gamma)}.$$

We assume below conditions that make this transformation well defined. By Itô's lemma, $\widetilde{X}_t \equiv g(X_t;\gamma)$ defined on $\widetilde{\mathcal{S}} = (g(\underline{x};\gamma), g(\bar{x};\gamma))$ satisfies $d\widetilde{X}_t =$

$\tilde{\mu}(\widetilde{X}_t; \beta) \, dt + dW_t$ with

$$\tilde{\mu}(x; \beta) \equiv \frac{\mu(g^{\text{inv}}(x; \gamma); \theta)}{\sigma(g^{\text{inv}}(x; \gamma); \gamma)} - \frac{1}{2} \frac{\partial \sigma(g^{\text{inv}}(x; \gamma); \gamma)}{\partial x},$$

where $g^{\text{inv}}$ denotes the reciprocal transformation. We also define the scale and speed densities of $\widetilde{X}$, $\tilde{s}$ and $\tilde{m}$, and $\tilde{\lambda}(x; \beta) \equiv -(\tilde{\mu}(x; \beta)^2 + \partial \tilde{\mu}(x; \beta)/\partial x)/2$.

We make the following primitive assumptions on $(\mu, \sigma)$:

ASSUMPTION 1. *For all values of the parameters $(\theta, \gamma)$ we have the following:*

1. *Differentiability*: The functions $\mu(x; \theta)$ and $\sigma(x; \gamma)$ are infinitely differentiable in $x$.
2. *Nondegeneracy of the diffusion*: If $\mathcal{S} = (-\infty, +\infty)$, there exists a constant $c$ such that $\sigma(x; \gamma) > c > 0$ for all $x$ and $\gamma$. If $\mathcal{S} = (0, +\infty)$, $\lim_{x \to 0^+} \sigma^2(x; \gamma) = 0$ is possible but then there exist constants $\xi_0 > 0$, $\omega > 0$, $\rho \geq 0$ such that $\sigma^2(x; \gamma) \geq \omega x^\rho$ for all $0 < x \leq \xi_0$ and $\gamma$. Whether or not $\lim_{x \to 0^+} \sigma^2(x; \gamma) = 0$, $\sigma$ is nondegenerate in the interior of $\mathcal{S}$, that is for each $\xi > 0$, there exists a constant $c_\xi$ such that $\sigma^2(x; \gamma) \geq c_\xi > 0$ for all $x \in [\xi, +\infty)$ and $\gamma$.
3. *Boundary behavior*: $\mu$, $\sigma^2$ and their derivatives have at most polynomial growth in $x$ near the boundaries, $\lim_{x \to \underline{x}} S(x; \beta) = -\infty$ and $\lim_{x \to \bar{x}} S(x; \beta) = +\infty$,

(8) $$\liminf_{x \to \underline{x}} \tilde{\mu}(x; \beta) > 0 \quad \text{and} \quad \limsup_{x \to \bar{x}} \tilde{\mu}(x; \beta) < 0$$

and

(9) $$\lim_{x \to \underline{x} \text{ or } x \to \bar{x}} \sup \tilde{\lambda}(x; \beta) < +\infty.$$

4. *Identification*: $\mu(x; \theta) = \mu(x; \tilde{\theta})$ for $\pi$-almost all $x$ in $\mathcal{S}$ implies $\theta = \tilde{\theta}$ and $\sigma^2(x; \gamma) = \sigma^2(x; \tilde{\gamma})$ for $\pi$-almost all $x$ in $\mathcal{S}$ implies $\gamma = \tilde{\gamma}$.

Under Assumption 1, the stochastic differential equation (1) admits a weak solution which is unique in probability law. This follows from the Engelbert–Schmidt criterion [see, e.g., Theorem 5.5.15 in Karatzas and Shreve (1991) replacing $\mathbb{R}$ by $\mathcal{S}$ throughout], with explosions ruled out by the boundary behavior of the process. The divergence of the scale measure makes the boundaries unattainable, since it implies that

$$\Sigma_{\bar{x}} \equiv \int^{\bar{x}} \left\{ \int_u^{\bar{x}} s(v; \beta) \, dv \right\} m(u; \beta) \, du = \infty,$$

$$\Sigma_{\underline{x}} \equiv \int_{\underline{x}} \left\{ \int_{\underline{x}}^u s(v; \beta) \, dv \right\} m(u; \beta) \, du = \infty.$$

Given that it is unattainable (i.e., given that $\Sigma_{\bar{x}} = \infty$), the boundary $\bar{x}$ is natural when $N_{\bar{x}} = \infty$ and entrance when $N_{\bar{x}} < \infty$, where

$$N_{\bar{x}} \equiv \int^{\bar{x}} \left\{ \int_u^{\bar{x}} m(v;\beta)\, dv \right\} s(u;\beta)\, du,$$

and similarly for the boundary $\underline{x}$ [see, e.g., Section 15.6 in Karlin and Taylor (1981)]. If both boundaries are entrance, then the integrability assumption on the speed measure $m$ will automatically be satisfied. When one of the boundaries is natural, integrability of $m$ is neither implied nor precluded.

Condition (8), however, guarantees that the process $\tilde{X}$ will be stationary: with $b$ denoting the lim sup in (8), we have near the right boundary

$$\tilde{m}(y;\beta) = \exp\left\{ 2 \int^x \tilde{\mu}(y;\theta)\, dy \right\} \leq c \exp\{2bx\},$$

where $c$ is a constant, and similarly near the left boundary. Thus, $\tilde{m}$ is integrable, and it follows that $m$ is integrable. Therefore, the process $X$ is stationary with stationary density

$$(10) \qquad \pi(x,\beta) = \frac{m(x;\beta)}{\int_{\underline{x}}^{\bar{x}} m(y;\beta)\, dy},$$

provided that the initial value of the process, $X_0$, has density $\pi$, which we will assume in the rest of the paper. Furthermore, condition (8) guarantees that $X$ has exponentially decaying $\rho$-mixing coefficients (see Lemma 4). Condition (9) guarantees the regularity of the transition density of the process [see Proposition 2 in Aït-Sahalia (2002)]; note that the condition does not prevent $\tilde{\lambda}$ from going to $-\infty$, it only excludes $+\infty$ as a possible limit.

We will denote by $L^2$ the Hilbert space of measurable real-valued functions $f$ on $\mathcal{S}$ such that $\|f\|^2 \equiv E[f(X_0)^2] < \infty$ for all values of $\beta$. When $f$ is a function of other variables, in addition to the state variable $y_1$, we say that $f \in L^2$ if it satisfies the integrability condition for every given value of the other variables.

To be able to give specific results on the effects of the sampling randomness on the estimation of $\beta$, we need to put some structure on the generation of the sampling intervals $\Delta_n = \tau_n - \tau_{n-1}$. We set $Y_n = X_{\tau_n}$. We assume the following regarding the data generating process for the sampling intervals:

ASSUMPTION 2. The sampling intervals $\Delta_n = \tau_n - \tau_{n-1}$ are independent and identically distributed. Each $\Delta_n$ is drawn from a common distribution which is independent of $Y_{n-1}$ and of the parameter $\beta$. Also, $E[\Delta_0^2] < +\infty$.

In particular, $E[f(Y_1)^2] = \|f\|^2$. An important special case occurs when the sampling happens to take place at a fixed deterministic interval $\bar{\Delta}$, corresponding to the distribution of $\Delta_n$ being a Dirac mass at $\bar{\Delta}$. See Section [5.3](#) for extensions. Throughout the paper we denote by $\Delta$ a generic random variable with the common distribution of the $\Delta_n$s.

While we assume that the distribution of the sampling intervals is independent of $\beta$, it may well depend upon its own nuisance parameters (such as an unknown arrival rate), but we are not interested in drawing inference about the (nuisance) parameters driving the sampling scheme, only about $\beta$. Note also that in the case of random sampling times, the number $N_T + 1$ of observations in the interval $[0, T]$ will be random.

2.2. *The estimators and their distribution.* We consider a class of estimators for $\beta$ obtained by minimizing a criterion function. Specifically, to estimate the $d$-dimensional parameter vector $\beta$, we select a vector of $r$ moment conditions $h(y_1, y_0, \delta, \beta, \varepsilon)$, $r \geq d$, which is continuously differentiable in $\beta$. We form the sample average

$$
(11) \qquad m_T(\beta) \equiv N_T^{-1} \sum_{n=1}^{N_T-1} h(Y_n, Y_{n-1}, \Delta_n, \beta, \varepsilon)
$$

and obtain $\hat{\beta}$ by minimizing the quadratic form

$$
(12) \qquad Q_T(\beta) \equiv m_T(\beta)' W_T m_T(\beta),
$$

where $W_T$ is an $r \times r$ positive definite weight matrix assumed to converge in probability to a positive definite limit $W_\beta$. If the system is exactly identified, $r = d$, the choice of $W_T$ is irrelevant and minimizing (12) amounts to setting $m_T(\beta)$ to 0. The function $h$ is known in different strands of the literature either as a "moment function" [see, e.g., Hansen (1982)] or an "estimating equation" [see, e.g., Godambe (1960) and Heyde (1997)].

A natural requirement on the estimating equation—what is needed for consistency of $\hat{\beta}$—is that

$$
(13) \qquad E_{\Delta, Y_1, Y_0}[h(Y_1, Y_0, \Delta, \beta_0, \varepsilon)] = 0.
$$

Throughout the paper we denote by $E_{\Delta, Y_1, Y_0}$ expectations taken with respect to the joint law of $(\Delta, Y_1, Y_0)$ at the true parameter $\beta_0$, and write $E_{\Delta, Y_1}$, and so on, for expectations taken from the appropriate marginal laws of $(\Delta, Y_1)$, and so on. As will become clear in the Euler example below, some otherwise fairly natural estimating strategies lead to inconsistent estimators. To allow for this, we do not assume that (13) is necessarily satisfied. Rather, we simply assume that the equation $E_{\Delta, Y_1, Y_0}[h(Y_1, Y_0, \Delta, \beta, \varepsilon)] = 0$ admits a unique root in $\beta$, which we define as $\bar{\beta} = \bar{\beta}(\beta_0, \varepsilon)$.

With $N_T/T$ converging in probability to $(E[\Delta])^{-1}$, it follows from standard arguments that $\sqrt{T}(\hat{\beta} - \bar{\beta})$ converges in law to $N(0, \Omega_\beta)$, with

$$
(14) \qquad \Omega_\beta^{-1} = (E[\Delta])^{-1} D_\beta' S_\beta^{-1} D_\beta,
$$

where

$$D_\beta \equiv E_{\Delta,Y_1,Y_0}[\dot{h}(Y_1, Y_0, \Delta, \bar{\beta}, \varepsilon)],$$
$$S_{\beta,j} \equiv E_{\Delta,Y_1,Y_0}[h(Y_{1+j}, Y_j, \Delta, \bar{\beta}, \varepsilon)h(Y_1, Y_0, \Delta, \bar{\beta}, \varepsilon)'],$$

and $S_\beta \equiv \sum_{j=-\infty}^{+\infty} S_{\beta,j}$. If $r > d$, the weight matrix $W_T$ is assumed to be any consistent estimator of $S_\beta^{-1}$; otherwise its choice is irrelevant. A consistent first-step estimator of $\bar{\beta}$, needed to compute the optimal weight matrix, can be obtained by minimizing (12) with $W_T = Id$.

## 3. Expansion of the asymptotic variance: general results.

3.1. *The generalized infinitesimal generator.* As we just saw, the asymptotic distributions of the estimators depend upon expectations of matrices of the form $E_{\Delta,Y_1,Y_0}[f(Y_1, Y_0, \Delta, \beta, \varepsilon)]$. Because these expectations are not generally available in closed-form, our approach is based on calculating Taylor expansions in $\varepsilon$ of these matrices. The key aspect of our approach is that these Taylor expansions all happen to be fully explicit.

To calculate Taylor expansions in $\varepsilon$ of the asymptotic variances when the sampling intervals are random, we introduce the *generalized infinitesimal operator* $\Gamma_{\beta_0}$ for the process $X$ in (1). This is in analogy to the development in Aït-Sahalia and Mykland (2003), but permits our current more general form of $\sigma^2(x;\gamma)$. To define this operator, let us first recall a standard concept. The standard infinitesimal generator $A_{\beta_0}$ is the operator which returns

$$(15) \quad A_{\beta_0} \cdot f = \frac{\partial f}{\partial \delta} + \mu(y_1, \theta_0)\frac{\partial f}{\partial y_1} + \frac{1}{2}\sigma^2(y_1;\gamma_0)\frac{\partial^2 f}{\partial y_1^2}$$

when applied to functions $f$ that are continuously differentiable once in $\delta$, twice in $y_1$ and such that $\partial f/\partial y_1$ and $A_{\beta_0} \cdot f$ are both in $L^2$ and satisfy

$$(16) \quad \lim_{y_1 \to \underline{x}} \frac{\partial f/\partial y_1}{s(y_1;\beta)} = \lim_{y_1 \to \bar{x}} \frac{\partial f/\partial y_1}{s(y_1;\beta)} = 0.$$

We define $\mathcal{D}$ to be the set of functions $f$ which have these properties and are additionally continuously differentiable in $\beta$ and $\varepsilon$. For instance, functions $f$ that are polynomial in $y_1$ near the boundaries of $\mathcal{S}$, and their iterates by repeated application of the generator, retain their polynomial growth characteristic near the boundaries; so they are all in $L^2$ and satisfy (16). Near both boundaries, polynomials and their iterates diverge at most polynomially (under Assumption 1, $\mu$, $\sigma^2$ and their derivatives have at most polynomial growth; multiplying and adding functions with polynomial growth yields a function still with polynomial growth). But we will often have exponential divergence of $s(y_1;\beta)$. This would be the case, for instance, if the left boundary is $\underline{x} = -\infty$, and there exist constants $E > 0$ and $K > 0$ such that for all $x < -E$ and $(\theta, \gamma)$, $\mu(x;\theta)/\sigma^2(x;\gamma) \geq K|x|^\alpha$ for some $\alpha \geq 0$; and if the right boundary is $\bar{x} = +\infty$, there exist constants $E > 0$ and $K > 0$ such that for all $x > E$ and $(\theta, \gamma)$, $\mu(x;\theta)/\sigma^2(x;\gamma) \leq -Kx^\alpha$ for some $\alpha \geq 0$; if instead $\underline{x} = 0^+$, there exist constants $E > 0$ and $K \geq 0$ such that for all $0 < x < E$ and $(\theta, \gamma)$, $\mu(x;\theta)/\sigma^2(x;\gamma) \geq Kx^{-\phi}$ for some $\phi > 1$ and $K > 0$. If, however, $\phi = 1$ and $K \geq 1$, then Assumption 1 is still satisfied, but $s$ diverges only polynomially near $0^+$.

Our new operator $\Gamma_{\beta_0}$ is then defined by its action on $f \in \mathcal{D}$:

$$(17) \qquad \Gamma_{\beta_0} \cdot f \equiv \Delta_0 A_{\beta_0} \cdot f + \frac{\partial f}{\partial \varepsilon} + \frac{\partial f}{\partial \beta} \frac{\partial \beta}{\partial \varepsilon}.$$

Note that when $\Delta_0$ is random, our operator $\Gamma_{\beta_0}$ is also random since it depends on $\Delta_0$. The last term allows for the fact that $\beta$ can be a function of $(\beta_0, \varepsilon)$. Because we will need to apply repeatedly the operator $\Gamma_{\beta_0}$, let us define $\mathcal{D}^J$ as the set of functions $f$ with $J+2$ continuous derivatives in $\delta$, $2(J+2)$ in $y_1$, such that $f$ and its first $J$ iterates by repeated applications of $A_{\beta_0}$ all remain in $\mathcal{D}$ and additionally have $J + 2$ continuous derivatives in $\beta$ and $\varepsilon$.

3.2. *Behavior for small $\varepsilon$ of the estimating equations.* The limiting behavior of the vector $h$ of moment functions depends crucially on whether one is estimating separately $\theta$, $\gamma$ or both together. If one is only estimating the drift parameters $\theta$, it will typically be the case that $h(y_1, y_0, \delta, \beta, \varepsilon)$ can be Taylor expanded around its continuous-sampling limit $h(y_0, y_0, 0, \beta, 0)$. On the other hand, when estimating $\gamma$, we will see that such a Taylor expansion is not possible, and $h(y_1, y_0, \delta, \beta, \varepsilon)$ is instead of order $O_p(1)$ as $\varepsilon \to 0$ (and naturally $y_1 \to y_0$ and $\delta \to 0$ at the same time). We shall in the following describe a structure which is consistent with all the various estimators of $\beta$ we consider, and is applicable to others as well.

We assume the following regularity condition regarding the moment functions $h$ selected to conduct inference:

ASSUMPTION 3. $h(y_1, y_0, \delta, \beta, \varepsilon) \in \mathcal{D}^J$ for some $J \geq 3$. We shall in general consider moment functions $h$ of the form

$$(18) \qquad h(y_1, y_0, \delta, \beta, \varepsilon) = \tilde{h}(y_1, y_0, \delta, \beta, \varepsilon) + \frac{H(y_1, y_0, \delta, \beta, \varepsilon)}{\delta},$$

where $\tilde{h} \in \mathcal{D}^J$ and $H \in \mathcal{D}^{J+1}$. When the function $H$ is not identically zero, we add the requirements that

$$(19) \quad H(y_0, y_0, 0, \beta_0, 0) = \frac{\partial H(y_1, y_0, 0, \beta_0, 0)}{\partial y_1} = \dot{H}(y_0, y_0, 0, \beta_0, 0) = 0.$$

In our definition of $h$, the term $H$ captures the singularity (i.e., powers of $1/\delta$) which can occur when estimating the diffusion coefficient. Consider, for example, the case where $h$ is the likelihood score for $\gamma$. The log-likelihood expansion in Aït-Sahalia (2002) is, for the transformed process $\tilde{X}$,

$$\tilde{l}(\tilde{y}_1, \tilde{y}_0, \delta, \beta, \varepsilon) = -\frac{1}{2}\ln(2\pi\delta) - \frac{(\tilde{y}_1 - \tilde{y}_0)^2}{2\delta} + \text{higher-order terms.}$$

By the Jacobian formula, the log-likelihood expansion for the original process $X$ is, therefore,

$$l(y_1, y_0, \delta, \beta, \varepsilon) = -\frac{1}{2}\ln(2\pi\delta) - \frac{1}{2\delta}\left(\int_{y_0}^{y_1} \frac{1}{\sigma(x;\gamma)}\,dx\right)^2$$
$$- \frac{1}{2}\ln(\sigma^2(y_1;\gamma)) + \text{higher-order terms.}$$

The score with respect to $\gamma$ is

$$h(y_1, y_0, \delta, \beta, \varepsilon) = \frac{1}{\delta}\left(\int_{y_0}^{y_1} \frac{1}{\sigma(x;\gamma)}\,dx\right)\left(\int_{y_0}^{y_1} \frac{\partial\sigma(x;\gamma)/\partial\gamma}{\sigma(x;\gamma)^2}\,dx\right)$$
$$- \frac{\partial\sigma(y_1;\gamma)/\partial\gamma}{\sigma(y_1;\gamma)} + \text{higher-order terms.}$$

$H$ must contain the coefficient of $\delta^{-1}$ in $h$, that is,

$$(20) \qquad H(y_1, y_0, \delta, \beta, \varepsilon) = \left(\int_{y_0}^{y_1} \frac{1}{\sigma(x;\gamma)}\,dx\right)\left(\int_{y_0}^{y_1} \frac{\partial\sigma(x;\gamma)/\partial\gamma}{\sigma(x;\gamma)^2}\,dx\right).$$

But we are free to add the terms of order $\delta$ and higher to $H$, provided we subtract them from $\tilde{h}$ so as to leave $h$ unchanged. For instance, a convenient choice is

$$(21) \quad \begin{aligned}&H(y_1, y_0, \delta, \beta, \varepsilon)\\ &= \left(\int_{y_0}^{y_1} \frac{1}{\sigma(x;\gamma)}\,dx\right)\left(\int_{y_0}^{y_1} \frac{\partial\sigma(x;\gamma)/\partial\gamma}{\sigma(x;\gamma)^2}\,dx\right) - \delta\frac{\partial\sigma(y_1;\gamma)/\partial\gamma}{\sigma(y_1;\gamma)}\end{aligned}$$

and then $\tilde{h} = h - H/\delta$. Both choices of $H$ satisfy (19) because $H$ has the form $H = a(y_1)b(y_1) - \delta c(y_1)$, where $a$ and $b$ denote, respectively, the two integrals and $c$ the coefficient of $\delta$ in (21); $c = 0$ in (20). At $y_1 = y_0$ and $\delta = 0$, we have $a = b = 0$, thus, $H = 0$; and $H' = a'b + ab' - \delta c' = 0$; and $\dot{H} = \dot{a}b + a\dot{b} - \delta\dot{c} = 0$. Note also that

$$A_{\beta_0} \cdot H = -c + \mu(a'b + ab' - \delta c') + (\sigma^2/2)(2a'b' + a''b + ab'' - \delta c''),$$

where prime denotes differentiation with respect to $y_1$. At $y_1 = y_0$ and $\delta = 0$, we have $a = b = 0$, thus, $A_{\beta_0} \cdot H = -c(y_0) + \sigma^2(y_0;\gamma_0)a'(y_0)b'(y_0) = 0$. And this $H$ does not depend on $\varepsilon$, so $\partial H/\partial\varepsilon = 0$, and the likelihood score is a martingale estimating function, hence unbiased, and so $\partial\bar{\beta}/\partial\varepsilon = 0$. Thus, by adding the additional term to $H$ in (21) relative to (20), we also have that $\Gamma_{\beta_0} \cdot H = 0$ at 0, which makes such an $H$ closer to being all by itself a martingale in a sense we make precise in Section 3.4. It makes no difference for the exact likelihood since $h$ is always a martingale, even if $\tilde{h}$ and $H$ are not separately, but this is convenient when analyzing likelihood approximations such as the Euler case in Section 4.2.

3.3. *The $D_\beta$ and $S_{\beta,0}$ matrices.* In order to obtain expansions of the form (3) for $\Omega_\beta$, we work on its components $D_\beta$, $S_{\beta,0}$ and $T_\beta = S_\beta - S_{\beta,0}$. The first result uses our generalized infinitesimal generator to provide the expansions of the matrices $D_\beta$ and $S_{\beta,0}$:

LEMMA 1 (Expansions for $D_\beta$ and $S_{\beta,0}$).   *Let $h = (h_1, \ldots, h_r)'$ denote a vector of moment functions $h = \tilde{h} + \Delta^{-1} H$ satisfying Assumption 3 and $\tilde{h} \in \mathcal{D}^3$, $H \in \mathcal{D}^4$. Also assume (27) and the other conditions on $q_i$ in Lemma 2.*

1. *In the case where $H$ is identically zero, we have*

$$(22) \qquad D_\beta = E_{Y_0}[\dot{h}] + \varepsilon E_{\Delta, Y_0}[(\Gamma_{\beta_0} \cdot \dot{h})] + \frac{\varepsilon^2}{2} E_{\Delta, Y_0}[(\Gamma_{\beta_0}^2 \cdot \dot{h})] + O(\varepsilon^3)$$

*and, with the notation $h \times h'(y_1, y_0, \delta, \beta, \varepsilon) \equiv h(y_1, y_0, \delta, \beta, \varepsilon) h(y_1, y_0, \delta, \beta, \varepsilon)'$, we have*

$$(23) \qquad \begin{aligned} S_{\beta,0} &= E_{Y_0}[(h \times h')] + \varepsilon E_{\Delta, Y_0}[(\Gamma_{\beta_0} \cdot (h \times h'))] \\ &\quad + \frac{\varepsilon^2}{2} E_{\Delta, Y_0}[(\Gamma_{\beta_0}^2 \cdot (h \times h'))] + O(\varepsilon^3). \end{aligned}$$

2. *In the case where $H$ is not zero, (22) and (23) should be evaluated at $\tilde{h}$ rather than $h$, yielding $D_\beta^{\tilde{h}}$ and $S_{\beta,0}^{\tilde{h}}$, respectively. Then $D_\beta = D_\beta^{\tilde{h}} + D_\beta^H$ and $S_{\beta,0} = S_{\beta,0}^{\tilde{h}} + S_{\beta,0}^H$, where*

$$(24) \quad D_\beta^H = E_{\Delta, Y_0}[\Delta_0^{-1}(\Gamma_{\beta_0} \cdot \dot{H})] + \frac{\varepsilon}{2} E_{\Delta, Y_0}[\Delta_0^{-1}(\Gamma_{\beta_0}^2 \cdot \dot{H})] + O(\varepsilon^2),$$

$$(25) \quad \begin{aligned} S_{\beta,0}^H &= E_{\Delta, Y_0}[\Delta_0^{-1}(\Gamma_{\beta_0} \cdot (\tilde{h} \times H'))] + \frac{\varepsilon}{2} E_{\Delta, Y_0}[\Delta_0^{-1}(\Gamma_{\beta_0}^2 \cdot (\tilde{h} \times H'))] \\ &\quad + E_{\Delta, Y_0}[\Delta_0^{-1}(\Gamma_{\beta_0} \cdot (H \times \tilde{h}'))] + \frac{\varepsilon}{2} E_{\Delta, Y_0}[\Delta_0^{-1}(\Gamma_{\beta_0}^2 \cdot (H \times \tilde{h}'))] \\ &\quad + \frac{1}{2} E_{\Delta, Y_0}[\Delta_0^{-2}(\Gamma_{\beta_0}^2 \cdot (H \times H'))] \\ &\quad + \frac{\varepsilon}{6} E_{\Delta, Y_0}[\Delta_0^{-2}(\Gamma_{\beta_0}^3 \cdot (H \times H'))] + O(\varepsilon^2). \end{aligned}$$

3.4. *How far is h from a martingale?* Next, we turn to an analysis of the more challenging time series matrix $T_\beta = S_\beta - S_{\beta,0}$. The simplest case arises when the moment function is a martingale,

$$(26) \qquad E_{\Delta,Y_1}[h(Y_1, Y_0, \Delta, \beta_0, \varepsilon)|Y_0] = 0.$$

In this circumstance, $S_{\beta,j} = 0$ for all $j \neq 0$, and so $T_\beta = 0$.

Even in the serially correlated case, however, we will show that the sum of these time series terms can, nonetheless, be small when $\varepsilon$ is small. Intuitively, the closer $h$ will be to a martingale, the smaller $T_\beta = S_\beta - S_{\beta,0}$.

To define what we mean by the distance from a moment function to a martingale, denote by $h_i$ the $i$th element of the vector of moment functions $h$, and define $q_i$ and $\alpha_i$ by

$$(27) \qquad \begin{aligned} & E_{\Delta,Y_1}[h_i(Y_1, Y_0, \Delta, \bar{\beta}, \varepsilon)|Y_0] \\ & \equiv \varepsilon^{\alpha_i} q_i(Y_0, \beta_0, \varepsilon) \\ & = \varepsilon^{\alpha_i} q_i(Y_0, \beta_0, 0) + \varepsilon^{\alpha_i+1} \frac{\partial q_i(Y_0, \beta_0, 0)}{\partial \varepsilon} + O(\varepsilon^{\alpha_i+2}), \end{aligned}$$

where $\alpha_i$ is an integer greater than or equal to zero for each moment function $h_i$. $\alpha_i$ is an index of the order at which the moment component $h_i$ deviates from a martingale (note that in a vector $h$ not all components $h_i$ need to have the same index $\alpha_i$). A martingale moment function corresponds to the limiting case where $\alpha_i = +\infty$, $q_i(Y_0, \beta_0, \varepsilon)$ is identically zero, and $S_\beta = S_{\beta,0}$. When the moment functions are not martingales, we will show that the difference $T_\beta \equiv S_\beta - S_{\beta,0}$ is a matrix whose element $(i,j)$ has a leading term of order $O(\varepsilon^{\min(\alpha_i,\alpha_j)})$ in $\varepsilon$ that depends on $q_i$ and $q_j$. As will become apparent in the following sections, (27) holds in all the estimation methods we consider.

Note that $E_{\Delta,Y_1,Y_0}[h_i(Y_1, Y_0, \Delta, \bar{\beta}, \varepsilon)] = 0$ by definition of $\bar{\beta}$, hence, by the law of iterated expectations we have that $E_{Y_0}[q_i(Y_0, \beta_0, \varepsilon)] = 0$. We will also need the function $r_i(y, \beta_0, \varepsilon)$ defined as

$$(28) \qquad r_i(y_0, \beta_0, \varepsilon) = -\int_0^\infty U_t \cdot A_{\beta_0} \cdot q_i(y_0, \beta_0, \varepsilon) E[\tau_{N(t)+1}] \, dt,$$

where $U_\delta \cdot f(y_0, \delta, \beta, \varepsilon) \equiv E_{Y_1}[f(Y_1, Y_0, \Delta, \beta, \varepsilon)|Y_0 = y_0, \Delta = \delta]$ is the conditional expectation operator. Recall that $\tau_i$ are the sampling times for $X$, and that $N_t = \#\{\tau_i \in (0, t]\}$ (so $\tau_0 = 0$ is not counted). We can assert the following about $r_i$:

LEMMA 2. *Under Assumption* 1 *and* (27), *we suppose that* $q_i(Y_0, \beta_0, 0)$ *and* $\frac{\partial q_i(Y_0,\beta_0,0)}{\partial \varepsilon}$ *are in* $L^2$. *Finally, let* $(A_{\beta_0} \cdot q_i)(Y, \beta_0, \varepsilon)$ *be defined, bounded and continuous in* $L^2$-*norm on an interval* $\varepsilon \in [0, \varepsilon_0]$ ($\varepsilon_0 > 0$). *Then* $r_i(y, \beta_0, \varepsilon)$ *is well defined. Also,*

$$(29) \qquad r_i(Y_0, \beta_0, \varepsilon) = \breve{r}_i(Y_0, \beta_0, \varepsilon) + \frac{1}{2}\varepsilon \frac{E[\Delta_0^2]}{E[\Delta_0]} q_i(Y_0, \beta_0, 0) + o_p(\varepsilon),$$

*where*

$$\begin{aligned}
\check{r}_i(y_0, \beta_0, \varepsilon) &= -\int_0^\infty t\ (U_t \cdot A_{\beta_0} \cdot q_i)(y_0, \beta_0, \varepsilon)\, dt \\
&= \int_0^\infty (U_t \cdot q_i)(y_0, \beta_0, \varepsilon)\, dt.
\end{aligned} \quad (30)$$

*Alternatively, $\check{r}_i$ can be defined as the solution of the differential equation*

$$A_{\beta_0} \cdot \check{r}_i(\cdot, \beta_0, \varepsilon) = -q_i(\cdot, \beta_0, \varepsilon), \quad (31)$$

*with the side condition that $E_{Y_0}[\check{r}_i(Y_0, \beta_0, \varepsilon)] = 0$.*

By convention, here and in the proofs $O_p(f(\varepsilon))$ and $o_p(f(\varepsilon))$ refer to terms whose $L^2$ norms are, respectively, $O(f(\varepsilon))$ and $o(f(\varepsilon))$.

Finally, an alternative form of $\check{r}_i$ is given in the proof of Lemma 2; see (64). While the index $\alpha_i$ and the function $q_i$ play a crucial role in determining the order in $\varepsilon$ of the matrix $T_\beta$, the function $r_i$ will play an important role in the determination of its coefficients.

3.5. *The $T_\beta$ matrix.* Putting all this together, we can calculate $T_\beta$ when $h$ is not a martingale estimating equation. The expansion of the matrix $T_\beta$ is obtained by applying the operator $\Gamma_{\beta_0}$ as follows:

LEMMA 3 (Expansions for $T_\beta$). *Under the assumptions of Lemma 1, assume also (27) and the other conditions on $q_i$ in Lemma 2. Then we have the following:*

1. *If $H$ is zero, the $(i,j)$ term of the time series matrix $T_\beta = S_\beta - S_{\beta,0}$ is $[T_\beta]_{(i,j)}$, given by*

$$\begin{aligned}
[T_\beta]_{(i,j)} = \frac{1}{E[\Delta_0]} \bigg( & \varepsilon^{\alpha_j - 1} E_{Y_0}[(h_i \times r_j)] + \varepsilon^{\alpha_j} E_{\Delta_0, Y_0}[\Gamma_{\beta_0} \cdot (h_i \times r_j)] \\
& + \frac{\varepsilon^{\alpha_j + 1}}{2} E_{\Delta_0, Y_0}[\Gamma_{\beta_0}^2 \cdot (h_i \times r_j)] + \varepsilon^{\alpha_i - 1} E_{Y_0}[(h_j \times r_i)] \\
& + \varepsilon^{\alpha_i} E_{\Delta_0, Y_0}[(\Gamma_{\beta_0} \cdot (h_j \times r_i))] \\
& \qquad\qquad + \frac{\varepsilon^{\alpha_i + 1}}{2} E_{\Delta_0, Y_0}[(\Gamma_{\beta_0}^2 \cdot (h_j \times r_i))] \bigg) \\
& + O(\varepsilon^{\min(\alpha_i, \alpha_j) + 2}).
\end{aligned} \quad (32)$$

2. If $H$ is nonzero, then (32) should be evaluated at $\tilde{h}$ rather than $h$, yielding $T_\beta^{\tilde{h}}$. And $T_\beta = T_\beta^{\tilde{h}} + T_\beta^H$, where

$$
\begin{aligned}
[T_\beta^H]_{(i,j)} = \frac{1}{E[\Delta_0]} &\bigg( \varepsilon^{\alpha_j - 1} E_{\Delta, Y_0}[\Delta_0^{-1}(\Gamma_{\beta_0} \cdot H_i) \times r_j] \\
&+ \frac{\varepsilon^{\alpha_j}}{2} E_{\Delta_0, Y_0}[\Delta_0^{-1}(\Gamma_{\beta_0}^2 \cdot (H_i \times r_j))] \\
&+ \varepsilon^{\alpha_i - 1} E_{\Delta_0, Y_0}[\Delta_0^{-1}(\Gamma_{\beta_0} \cdot H_j) \times r_i] \\
&+ \frac{\varepsilon^{\alpha_i}}{2} E_{\Delta_0, Y_0}[\Delta_0^{-1}(\Gamma_{\beta_0}^2 \cdot (H_j \times r_i))] \bigg) \\
&+ O(\varepsilon^{\min(\alpha_i, \alpha_j) + 1}).
\end{aligned}
\tag{33}
$$

Note that for most applications of Lemmas 1–3, the assumption of (27) and the other conditions on $q_i$ in Lemma 2 follow from the other assumptions of Lemma 3, as follows. Normally, one can take $\alpha_i$ to be $\le 2$, since the error term in (32) need not be smaller than that of (23), and the error term in (33) need not be smaller than that of (25). The conditions mentioned from Lemma 2 follow if $\tilde{h}_i \in \mathcal{D}^3$ and $H_i \in \mathcal{D}^4$ (more generally, $\tilde{h}_i \in \mathcal{D}^{\alpha_i + 1}$ and $H_i \in \mathcal{D}^{\alpha_i + 2}$).

3.6. *Form of the asymptotic variance matrix $\Omega_\beta$.* By combining our previous results concerning $D_\beta$, $S_{\beta,0}$ and $T_\beta$, we can now obtain an expression for the matrix $\Omega_\beta$. Specifically, we have the following as an example of a typical situation:

THEOREM 1 (Form of the matrix $\Omega_\beta$). *Under the conditions of the preceding lemmas we have the following:*

1. *When we are only estimating $\theta$, with $\gamma_0$ known, using a vector $h$ such that $H = 0$, and $D_\beta$ and $S_\beta$ have the expansions*

$$D_\theta = \varepsilon D_\theta^{(1)} + \varepsilon^2 D_\theta^{(2)} + O(\varepsilon^3),$$

$$S_\theta = \varepsilon S_\theta^{(1)} + \varepsilon^2 D_\theta^{(2)} + O(\varepsilon^3),$$

*then the asymptotic variance of the estimator has the expansion $\Omega_\theta = \Omega_\theta^{(0)} + \varepsilon \Omega_\theta^{(1)} + O(\varepsilon^2)$, where*

$$\Omega_\theta^{(0)} = E[\Delta_0] S_\theta^{(1)} / (D_\theta^{(1)})^2,$$

$$\Omega_\theta^{(1)} = E[\Delta_0](D_\theta^{(1)} S_\theta^{(2)} - 2 D_\theta^{(2)} S_\theta^{(1)}) / (D_\theta^{(1)})^3.$$

2. *When we are only estimating $\gamma$, with $\theta_0$ known, and $D_\beta$ and $S_\beta$ have the expansions*

$$D_\gamma = D_\gamma^{(0)} + \varepsilon D_\gamma^{(1)} + \varepsilon^2 D_\gamma^{(2)} + O(\varepsilon^3),$$

$$S_\gamma = S_\gamma^{(0)} + \varepsilon S_\gamma^{(1)} + \varepsilon^2 S_\gamma^{(2)} + O(\varepsilon^3),$$

*then the asymptotic variance of the estimator has the expansion* $\Omega_\gamma = \Omega_\gamma^{(1)}\varepsilon + \Omega_\gamma^{(2)}\varepsilon^2 + \Omega_\gamma^{(3)}\varepsilon^3 + O(\varepsilon^4)$, *where*

$$\Omega_\gamma^{(1)} = E[\Delta_0]S_\gamma^{(0)}/(D_\gamma^{(0)})^2,$$

$$\Omega_\gamma^{(2)} = E[\Delta_0](D_\gamma^{(0)}S_\gamma^{(1)} - 2D_\gamma^{(1)}S_\gamma^{(0)})/(D_\gamma^{(0)})^3,$$

$$\Omega_\gamma^{(3)} = E[\Delta_0](3D_\gamma^{(1)}S_\gamma^{(0)} - 2D_\gamma^{(0)}D_\gamma^{(2)}S_\gamma^{(2)}$$
$$- 2D_\gamma^{(0)}D_\gamma^{(1)}S_\gamma^{(1)} + (D_\gamma^{(0)})^2 S_\gamma^{(2)})/(D_\gamma^{(0)})^4.$$

3. *When we are estimating* $\theta$ *and* $\gamma$ *jointly, and the* $D_\beta$ *and* $S_\beta$ *matrices have the expansions*

$$D_\beta = \begin{pmatrix} \varepsilon d_{\theta\theta}^{(1)} + \varepsilon^2 d_{\theta\theta}^{(2)} + O(\varepsilon^3) & \varepsilon d_{\theta\gamma}^{(1)} + O(\varepsilon^2) \\ \varepsilon d_{\gamma\theta}^{(1)} + O(\varepsilon^2) & d_{\gamma\gamma}^{(0)} + \varepsilon d_{\gamma\gamma}^{(1)} + O(\varepsilon^2) \end{pmatrix},$$

$$S_\beta = \begin{pmatrix} \varepsilon s_{\theta\theta}^{(1)} + \varepsilon^2 s_{\theta\theta}^{(2)} + O(\varepsilon^3) & \varepsilon s_{\theta\gamma}^{(1)} + O(\varepsilon^2) \\ \varepsilon s_{\gamma\theta}^{(1)} + O(\varepsilon^2) & s_{\gamma\gamma}^{(0)} + \varepsilon s_{\gamma\gamma}^{(1)} + O(\varepsilon^2) \end{pmatrix},$$

*then the asymptotic variance of the estimator* $\hat{\beta}$ *is*

$$\Omega_\beta = \begin{pmatrix} \omega_{\theta\theta}^{(0)} + \varepsilon\omega_{\theta\theta}^{(1)} + O(\varepsilon^2) & \varepsilon\omega_{\theta\gamma}^{(1)} + O(\varepsilon^2) \\ \varepsilon\omega_{\gamma\theta}^{(1)} + O(\varepsilon^2) & \varepsilon\omega_{\gamma\gamma}^{(1)} + \varepsilon^2\omega_{\gamma\gamma}^{(2)} + O(\varepsilon^3) \end{pmatrix},$$

*where*

$$\omega_{\theta\theta}^{(0)} = E[\Delta_0]s_{\theta\theta}^{(1)}/(d_{\theta\theta}^{(1)})^2,$$

$$\omega_{\theta\theta}^{(1)} = E[\Delta_0](d_{\theta\theta}^{(1)}s_{\theta\theta}^{(2)} - 2d_{\theta\theta}^{(2)}s_{\theta\theta}^{(1)})/(d_{\theta\theta}^{(1)})^3,$$

$$\omega_{\gamma\theta}^{(1)} = \omega_{\theta\gamma}^{(1)} = E[\Delta_0]s_{\theta\gamma}^{(1)}/(d_{\theta\theta}^{(1)}d_{\gamma\gamma}^{(0)}),$$

$$\omega_{\gamma\gamma}^{(1)} = E[\Delta_0]s_{\gamma\gamma}^{(0)}/(d_{\gamma\gamma}^{(0)})^2,$$

$$\omega_{\gamma\gamma}^{(2)} = E[\Delta_0](d_{\gamma\gamma}^{(0)}d_{\gamma\gamma}^{(1)} - 2d_{\gamma\gamma}^{(1)}s_{\gamma\gamma}^{(0)})/(d_{\gamma\gamma}^{(0)})^3.$$

*In particular, the diagonal leading terms* $\omega_{\theta\theta}^{(0)}$ *and* $\omega_{\gamma\gamma}^{(1)}$ *corresponding to efficient estimation with a continuous record of observation are identical to their single-parameter counterparts.*

An important fact to note from the above expressions is that to first order in $\varepsilon$, the asymptotic variances of $\hat{\theta}$ and $\hat{\gamma}$ are unaffected by whether one estimates just one of them (and the other one is known) or one estimates both of them jointly. This is not necessarily the case for the higher-order terms in the asymptotic variances, since those depend upon the higher-order terms in the $D_\beta$ and $S_\beta$ matrices which are not necessarily identical to those of their single-parameter counterparts.

Also, the leading term in $\Omega_\theta$ corresponding to efficient estimation with a continuous record of observations is

$$\Omega_\theta^{(0)} = (E_{Y_0}[(\partial \mu(Y_0;\theta_0)/\partial \theta)^2 \sigma(Y_0;\gamma_0)^{-2}])^{-1}$$

provided $\mu$ is continuously differentiable with respect to $\theta$. And the leading term in $\Omega_\gamma$ corresponding to efficient estimation of $\gamma$ is

$$\Omega_\gamma^{(1)} = E[\Delta_0](2 E_{Y_0}[(\partial \sigma(Y_0;\gamma_0)/\partial \gamma)^2 \sigma(Y_0;\gamma_0)^{-2}])^{-1}$$

provided $\sigma$ is continuously differentiable with respect to $\gamma$. In the special case where $\sigma^2 = \gamma$ constant, then this becomes $\Omega_{\sigma^2}^{(1)} = 2\sigma_0^4 E[\Delta_0]$.

These leading terms are achieved, in particular, when $h$ is the likelihood score for $\theta$ and $\gamma$, respectively, but also by other estimating functions that are able to mimic the behavior of the likelihood score at the leading order.

3.7. *Inconsistency.* For the estimator to be consistent, it must be that $\bar{\beta} \equiv \beta_0$ but, again, this will not be the case for every estimation method. However, in all the cases we consider, and one may argue for *any* reasonable estimation method, the bias will disappear in the limit where $\varepsilon \to 0$, that is, $\bar{\beta}(\beta_0, 0) = \beta_0$ (so that there is no bias in the limiting case of continuous sampling) and the following expansion

$$(34) \qquad \bar{\beta} = \bar{\beta}(\beta_0, \varepsilon) = \beta_0 + \sum_{q=1}^{Q} \varepsilon^q b^{(q)} + o(\varepsilon^Q)$$

holds for some $Q \geq 1$. The coefficients $b^{(q)} = (1/q!) \partial^q \bar{\beta}(\beta_0, 0)/\partial \varepsilon^q$ can be determined as follows. By the definition of $\bar{\beta}$,

$$(35) \qquad E_{\Delta, Y_1, Y_0}[h(Y_1, Y_0, \Delta, \bar{\beta}(\beta_0, \varepsilon), \varepsilon)] \equiv 0.$$

Consider the case where $H = 0$. Recognizing that $\bar{\beta}$ is a function of $\varepsilon$, as given in (34), we can compute the Taylor series expansion

$$(36) \qquad \begin{aligned} & E_{Y_1}[h(Y_1, Y_0, \Delta, \bar{\beta}, \varepsilon)|Y_0, \Delta] \\ & = \sum_{j=0}^{J} \frac{\varepsilon^j}{j!} (\Gamma_{\beta_0}^j \cdot h)(Y_0, Y_0, 0, \beta_0, 0) + O_p(\varepsilon^{J+1}), \end{aligned}$$

whose unconditional expectation, in light of (35), must be zero at each order in $\varepsilon$. So to determine $b^{(1)}$, set to zero the coefficient of $\varepsilon$ in the series expansion of $E_{\Delta, Y_1, Y_0}[h(Y_1, Y_0, \Delta, \bar{\beta}(\beta_0, \varepsilon), \varepsilon)] = E_{\Delta, Y_0}[E_{Y_1}[h(Y_1, Y_0, \Delta, \bar{\beta}(\beta_0, \varepsilon), \varepsilon)|Y_0, \Delta]]$:

$$\begin{aligned} 0 &= E_{\Delta, Y_0}[(\Gamma_{\beta_0} \cdot h)(Y_0, Y_0, 0, \beta_0, 0)] \\ &= E[\Delta_0] E_{Y_0}[(A_{\beta_0} \cdot h)(Y_0, Y_0, 0, \beta_0, 0)] \\ &\quad + E_{Y_0}\left[\frac{\partial h}{\partial \varepsilon}(Y_0, Y_0, 0, \beta_0, 0)\right] + E_{Y_0}[\dot{h}(Y_0, Y_0, 0, \beta_0, 0)] b^{(1)} \end{aligned}$$

and, hence, if $E_{Y_0}[\dot{h}(Y_0, Y_0, 0, \beta_0, 0)] \neq 0$,

$$(37) \quad \begin{aligned} b^{(1)} = -(&E[\Delta_0] E_{Y_0}[(A_{\beta_0} \cdot h)(Y_0, Y_0, 0, \beta_0, 0)] \\ &+ E_{Y_0}[(\partial h/\partial \varepsilon)(Y_0, Y_0, 0, \beta_0, 0)])(E_{Y_0}[\dot{h}(Y_0, Y_0, 0, \beta_0, 0)])^{-1}. \end{aligned}$$

Then given $b^{(1)}$, setting the coefficient of $\varepsilon^2$ in that series expansion to zero determines $b^{(2)}$, and so on. If $E_{Y_0}[\dot{h}(Y_0, Y_0, 0, \beta_0, 0)] = 0$, then one needs to look at the next order term in the expansion to determine $b^{(1)}$, and so on. This is, for instance, what happens in the Euler scheme when estimating $\theta$; see Section 4.2.

If $H \neq 0$, then (36) incorporates both $\tilde{h}$ and $H$, and one proceeds analogously to determine $b^{(1)}$ and the following coefficients by setting the coefficients of the expansion of (35) to 0. For an example of this, see the estimation of $\sigma^2$ using the Euler scheme.

**4. Application to specific inference strategies.** We now apply the general results to specific instances of moment functions $h$, corresponding both to likelihood and nonlikelihood inference strategies, for the model where $\sigma^2 = \gamma$ constant.

4.1. *Maximum-likelihood type estimators.* The development of Aït-Sahalia and Mykland (2003) deals with likelihood type inference, and we recapitulate here the inference schemes in that work, and how they relate to the present paper. We applied the general results of the present paper to maximum likelihood estimation, using three different inference strategies:

1. FIML: *Full information* maximum likelihood, using the bivariate observations $(Y_n, \Delta_n)$.
2. IOML: Partial information maximum likelihood estimator using only the state observations $Y_n$, with the sampling intervals *integrated out*.
3. PFML: Pseudo maximum likelihood estimator *pretending* that the sampling intervals are *fixed* at $\Delta_n = \bar{\Delta}$.

All three estimators rely on maximizing a version of the likelihood function of the observations, that is, some functional of the transition density $p$ of the $X$ process: $p(Y_n|Y_{n-1}, \Delta_n, \theta)$ for FIML; $\tilde{p}(Y_n|Y_{n-1}, \theta) = E_{\Delta_n}[p(Y_n|Y_{n-1}, \Delta_n, \theta]$ for IOML; and $p(Y_n|Y_{n-1}, \bar{\Delta}, \theta)$ for PFML (which is like FIML except that $\bar{\Delta}$ is used in place of the actual $\Delta_n$). The extent to which these estimators differ from one another gave rise to different "costs." FIML is asymptotically efficient, making the best possible use of the joint discretely sampled data $(Y_n, \Delta_n)$. The extent to which FIML with these data is less efficient than the corresponding FIML when the full sample path is observable is what we called the *cost of discreteness*. IOML is the asymptotically optimal choice if one recognizes that the sampling intervals $\Delta_n$ are random but does not observe them. The extra efficiency loss relative to FIML is what we called the *cost of randomness*. PFML corresponds to the "head-in-the-sand" policy consisting of doing as if the sampling intervals were all identical (pretending that $\Delta_n = \bar{\Delta}$) when, in fact, they are random. The extent by which PFML underperforms FIML is what we called the *cost of ignoring the randomness*. We then studied the relative magnitude of these costs in various situations.

The respective scores from these likelihoods are special cases of the estimating functions $h$ of the present paper. But the results of the present paper apply to a much wider class of estimating functions than the three likelihood examples, such as the following.

4.2. *Estimator based on the discrete Euler scheme.* We now apply our general results to study the properties of estimators of the drift and diffusion coefficients obtained by replacing the true likelihood function $l(y_1|y_0, \delta, \beta)$ with its discrete Euler approximation

$$(38) \qquad l^E(y_1|y_0, \delta, \beta) = -\frac{1}{2}\ln(2\pi\sigma^2\delta) - \frac{(y_1 - y_0 - \mu(y_0;\theta)\delta)^2}{2\sigma^2\delta}.$$

This estimator is commonly used in empirical work in finance, where researchers often write a theoretical model set in continuous-time but then switch gear in their empirical work, in effect estimating the parameters $\bar{\beta} = (\bar{\theta}, \bar{\sigma}^2)'$ of the discrete time series model

$$(39) \qquad X_{t+\Delta} - X_t = \mu(X_t; \bar{\theta})\Delta + \bar{\sigma}\sqrt{\Delta}\eta_{t+\Delta},$$

where the disturbance $\eta$ is $N(0,1)$. The properties of this estimator have been studied in the case where $\Delta$ is not random by Florens-Zmirou (1989). Our results apply to this particular situation as a special case.

In the terminology of Section 3, our vector of moment functions is

$$(40) \qquad \begin{aligned} h(y_1, y_0, \delta, \beta, \varepsilon) &= \begin{bmatrix} \dot{l}^E_\theta(y_1|y_0, \delta, \beta) \\ \dot{l}^E_{\sigma^2}(y_1|y_0, \delta, \beta) \end{bmatrix} \\ &= \begin{bmatrix} \dot{\mu}(y_0, \theta)(y_1 - y_0 - \mu(y_0;\theta)\delta)/\sigma^2 \\ -1/(2\sigma^2) + (y_1 - y_0 - \mu(y_0;\theta)\delta)^2/(2\sigma^4\delta) \end{bmatrix} \end{aligned}$$

when both parameters in $\beta = (\theta, \sigma^2)$ are unknown, and reduces to one component when only one parameter is unknown. For this choice of $h$, (13) is not satisfied and, thus, the estimator is inconsistent. Note also that the solution in $\theta$ of $E_{\Delta, Y_1, Y_0}[\dot{l}^E_\theta(Y_1, Y_0, \Delta, \beta, \varepsilon)] = 0$ is independent of $\sigma^2$ and, hence, whether or not we are estimating $\sigma^2$ does not affect the estimator of the drift parameter. Of course, this will not be the case in general for the true maximum likelihood estimator.

As we discussed in the general case, the asymptotic bias of the estimator, $\bar{\beta} - \beta_0$, will be of order $O(\varepsilon)$ or smaller. In this particular case, if $\sigma^2$ is known, the bias in $\hat{\theta}$ is of order $O(\varepsilon)$. As in the general setting of Section 3, $\sqrt{T}(\hat{\beta} - \bar{\beta})$ converges in law to $N(0, \Omega_\beta)$ and an application of Lemmas 1 and 3 yields the following.

1. When we are only estimating $\theta$, with $\sigma_0^2$ known, using only the first equation in (40), we have $\alpha_1 = 2$ and

$$
\begin{aligned}
&q_1(y, \beta_0, 0) \\
&= \frac{E[\Delta_0^2]}{4\sigma_0^2} \\
&\quad \times \left( \frac{\sigma_0^2 E_{Y_0}[(\partial\mu/\partial y)(Y_0; \theta_0)\,(\partial^2\mu/(\partial y\,\partial\theta))(Y_0; \theta_0)]((\partial\mu/\partial\theta)(y; \theta_0))^2}{E_{Y_0}[((\partial\mu/\partial\theta)(Y_0; \theta_0))^2]} \right. \\
&\quad\quad\left. + \left(2\mu(y; \theta_0)\frac{\partial\mu(y; \theta_0)}{\partial y} + \sigma_0^2 \frac{\partial^2\mu(y; \theta_0)}{\partial y^2}\right)\frac{\partial\mu(y; \theta_0)}{\partial\theta} \right)
\end{aligned}
\tag{41}
$$

in (27). The bias of the drift estimator is

$$
\begin{aligned}
&\bar{\theta} - \theta_0 \\
&= -\varepsilon\sigma_0^2 \frac{E[\Delta_0^2]}{E[\Delta_0]} \frac{E_{Y_0}[(\partial\mu/\partial y)(Y_0; \theta_0)(\partial^2\mu/(\partial y\,\partial\theta))(Y_0; \theta_0)]}{4 E_{Y_0}[((\partial\mu/\partial\theta)(Y_0; \theta_0))^2]} + O(\varepsilon^2)
\end{aligned}
\tag{42}
$$

and its asymptotic variance is $\Omega_\theta = \Omega_\theta^{(0)} + \Omega_\theta^{(1)}\varepsilon + O(\varepsilon^2)$ with $\Omega_\theta^{(0)} = \sigma_0^2 E_{Y_0}[((\partial\mu/\partial\theta)(Y_0; \theta_0))^2]^{-1}$ (the limiting term corresponding to a continuous record of observations) and

$$
\begin{aligned}
\Omega_\theta^{(1)} &= \frac{\sigma_0^2 E[\Delta_0^2]}{2 E[\Delta_0] E_{Y_0}[((\partial\mu/\partial\theta)(Y_0; \theta_0))^2]^3} \\
&\quad \times \left( 2\sigma_0^2 E_{Y_0}\left[\frac{\partial\mu}{\partial\theta}(Y_0; \theta_0)\frac{\partial^2\mu}{\partial\theta^2}(Y_0; \theta_0)\right] E_{Y_0}\left[\frac{\partial\mu}{\partial y}(Y_0; \theta_0)\frac{\partial^2\mu}{\partial y\,\partial\theta}(Y_0; \theta_0)\right] \right. \\
&\quad\quad + E_{Y_0}\left[\left(\frac{\partial\mu}{\partial\theta}(Y_0; \theta_0)\right)^2\right] \\
&\quad\quad \times \left( \frac{2\sigma_0^2 T_\theta^{(2)}}{E[\Delta_0^2]} + 2 E_{Y_0}\left[\left(\frac{\partial\mu}{\partial\theta}(Y_0; \theta_0)\right)^2 \frac{\partial\mu}{\partial y}(Y_0; \theta_0)\right] \right. \\
&\quad\quad\quad\left.\left.- \sigma_0^2 E_{Y_0}\left[\frac{\partial\mu}{\partial y}(Y_0; \theta_0)\frac{\partial^3\mu}{\partial y\,\partial\theta^2}(Y_0; \theta_0)\right]\right)\right),
\end{aligned}
$$

where

$$
T_\theta^{(2)} = -2 E_{Y_0}[\dot{q}_1(Y_0, \beta_0, 0)] = 4 E_{Y_0}[q_1(Y_0, \beta_0, 0) G_1(Y_0, \beta_0)]
$$

with $G_1(y_0, \beta_0) = \sigma_0^{-2} \int^{y_0} \dot{\mu}(z_0, \theta_0)\,dz_0$.

2. When we are only estimating $\sigma^2$, with $\theta_0$ known, using only the second equation in (40), we have $\alpha_2 = 1$ and

$$
q_2(y, \beta_0, 0) = \frac{E[\Delta_0]}{2\sigma_0^2}\left(\frac{\partial\mu}{\partial y}(y, \theta_0) - E_{Y_0}\left[\frac{\partial\mu}{\partial y}(Y_0; \theta_0)\right]\right)
\tag{43}
$$

in (27). The bias of the diffusion estimator is

$$\bar{\sigma}^2 - \sigma_0^2 = \varepsilon E[\Delta_0]\sigma_0^2 E_{Y_0}\left[\frac{\partial \mu}{\partial y}(Y_0;\theta_0)\right]$$
$$(44)\qquad + \varepsilon^2 \frac{2\sigma_0^2}{3} E[\Delta_0^2] E_{Y_0}\left[\left(\frac{\partial \mu}{\partial y}(Y_0;\theta_0)\right)^2\right] + O(\varepsilon^3)$$

and its asymptotic variance is $\Omega_{\sigma^2} = \Omega_{\sigma^2}^{(1)}\varepsilon + \Omega_{\sigma^2}^{(2)}\varepsilon^2 + O(\varepsilon^3)$ with $\Omega_{\sigma^2}^{(1)} = 2\sigma_0^4 E[\Delta_0]$ (the same first-order term as MLE) and

$$\Omega_{\sigma^2}^{(2)} = 4\sigma_0^4 E[\Delta_0]\left(E[\Delta_0]E_{Y_0}\left[\frac{\partial \mu}{\partial y}(Y_0;\theta_0)\right] + \sigma_0^4 T_{\sigma^2}^{(1)}\right),$$

where

$$T_{\sigma^2}^{(1)} = 4E_{Y_0}[q_2(Y_0,\beta_0,0)G_2(Y_0,\beta_0)] + \frac{2}{E[\Delta_0]}E_{Y_0}[q_2(Y_0,\beta_0,0)r_2(Y_0,\beta_0,0)]$$

with $G_2(y_0,\beta_0) = -\sigma_0^{-4}\int^{y_0} \mu(z_0,\theta_0)\,dz_0$.

3. When we are estimating $\theta$ and $\sigma^2$ jointly, using both equations in (40), the two components of the bias vector $\bar{\beta} - \beta_0$ are given by (42) and (44), respectively, (to their respective orders only). We also have that $\alpha_1 = 2, \alpha_2 = 1$ and $q = (q_1, q_2)'$ with $q_1$ and $q_2$ given by (41) and (43), respectively. The asymptotic variance of $\hat{\beta}$ is

$$\Omega_\beta = \begin{pmatrix} \omega_{\theta\theta} & \omega_{\theta\sigma^2} \\ \omega_{\sigma^2\theta} & \omega_{\sigma^2\sigma^2} \end{pmatrix}$$
$$= \begin{pmatrix} \omega_{\theta\theta}^{(0)} + \varepsilon\omega_{\theta\theta}^{(1)} + O(\varepsilon^2) & \varepsilon\omega_{\theta\sigma^2}^{(1)} + O(\varepsilon^2) \\ \varepsilon\omega_{\sigma^2\theta}^{(1)} + O(\varepsilon^2) & \varepsilon\omega_{\sigma^2\sigma^2}^{(1)} + \varepsilon^2\omega_{\sigma^2\sigma^2}^{(2)} + O(\varepsilon^3) \end{pmatrix},$$

where $\omega_{\theta\theta}^{(0)} = \Omega_\theta^{(0)}$, $\omega_{\theta\theta}^{(1)} = \Omega_\theta^{(1)}$, $\omega_{\sigma^2\sigma^2}^{(1)} = \Omega_{\sigma^2}^{(1)}$ and

$$\omega_{\sigma^2\theta}^{(1)} = \omega_{\theta\sigma^2}^{(1)} = \frac{2\sigma_0^6}{E_{Y_0}[(\partial\mu/\partial\theta(Y_0;\theta_0))^2]} t_{\theta\sigma^2}^{(1)},$$

$$\omega_{\sigma^2\sigma^2}^{(2)} = 4\sigma_0^4 E[\Delta_0]\left(E[\Delta_0]E_{Y_0}\left[\frac{\partial\mu}{\partial y}(Y_0;\theta_0)\right] + \sigma_0^4 t_{\sigma^2\sigma^2}^{(1)}\right),$$

with $t_{\theta\sigma^2}^{(1)} = 2E_{Y_0}[G_1(Y_0,\beta_0)q_2(Y_0,\beta_0,0)]$ and $t_{\sigma^2\sigma^2}^{(1)} = T_{\sigma^2}^{(1)}$.

Therefore, as is to be expected when using a first-order approximation to the stochastic differential equation, the asymptotic variance is, to first order in $\varepsilon$, the same as for MLE inference. The impact of using the approximation is to second order in variances (and, of course, is responsible for bias in the estimator). When estimating one of the two parameters with the other known, the impact of the discretization approximation on the variance (which MLE avoids) is one order of magnitude higher than the effect of the discreteness of the data (which MLE is also subject to).

4.3. *Example*: *the Ornstein–Uhlenbeck process.* We now specialize the expressions above to a specific example, the stationary ($\theta > 0$) Ornstein–Uhlenbeck process

$$dX_t = -\theta X_t\, dt + \sigma\, dW_t. \tag{45}$$

The transition density $p(y_1|y_0, \delta, \beta)$ of this process is a Gaussian density with expected value $e^{-\delta\theta}y_0$ and variance $(1 - e^{-2\delta\theta})\sigma^2/(2\theta)$. The stationary density $\pi(y_0, \beta)$ is also Gaussian with mean 0 and variance $\sigma^2/(2\theta)$.

Because its transition density is known explicitly, this model constitutes one of the rare instances where, in addition to our Taylor expansions which can be calculated for any model, we can obtain exact (i.e., non-Taylor expanded) expressions for the matrices $S_{\beta,0}$, $D_\beta$ and $T_\beta$. Specifically, for methods relying on nonmartingale moment functions $h$, the exact calculation of the time series term $T_\beta$ is based on

$$T_\beta = \frac{2}{E[\Delta_0]} E_{\Delta, Y_1, Y_0}[h(Y_1, Y_0, \Delta, \bar{\beta}, \varepsilon) R(Y_1, \beta_0, \varepsilon)],$$

where $E_{\Delta, Y_1}[h(Y_0, Y_1, \Delta, \bar{\beta}, \varepsilon)|Y_0] = \varepsilon^\alpha q(Y_0, \beta_0, \varepsilon) \equiv Q(Y_0, \beta_0, \varepsilon)$ and

$$R(Y_1, \beta_0, \varepsilon) = E[\Delta_0] \sum_{k=1}^\infty E_{Y_k}[Q(Y_k, \beta_0, \varepsilon)|Y_1] = \varepsilon^{\alpha-1} r(Y_1, \beta_0, \varepsilon).$$

This last expression requires the calculation of $E[Y_k^2|Y_1]$. To this end, consider first the law of $Y_k$ given $Y_1$ and $\Delta_2, \ldots, \Delta_k$. In this case, $Y_k$ is conditionally Gaussian with mean $Y_1 \exp\{-\theta(\Delta_2 + \cdots + \Delta_k)\}$ and variance $((k-1) - \exp\{-2\theta(\Delta_2 + \cdots + \Delta_k)\})\sigma^2/(2\theta)$. Hence, we obtain that

$$\begin{aligned}
E[Y_k^2|Y_1] &= E\bigg[Y_1^2 \exp\{-2\theta(\Delta_2 + \cdots + \Delta_k)\} \\
&\quad + \frac{\sigma^2}{2\theta}((k-1) - \exp\{-2\theta(\Delta_2 + \cdots + \Delta_k)\})|Y_1\bigg] \\
&= Y_1^2 E[\exp\{-2\theta\Delta\}]^{(k-1)} + \frac{\sigma^2}{2\theta}((k-1) - E[\exp\{-2\theta\Delta\}]^{(k-1)}).
\end{aligned}$$

In Table 1, we report results for the Ornstein–Uhlenbeck parameters estimated one at a time (i.e., $\theta$ knowing $\sigma^2$ and $\sigma^2$ knowing $\theta$). The quantities for the MLE are based on the developments in Aït-Sahalia and Mykland (2003); for the discrete Euler scheme, they follow from the results above.

4.4. *The effect of the distribution of the sampling intervals.* One of the implications of our results concerns the impact of the distribution of the sampling interval on the quality of inference. It is, obviously, better to have as many sampling times as possible, but, to move beyond this, fix $E[\Delta_0]$.

To the extent that our expansions depend on other features of the law of $\Delta_0$, they do so through the moments $E[\Delta_0^q]$, $q \geq 2$, as can be seen from the expressions above.

One can then compare whether it seems preferable to minimize these higher-order moments, and thus have sampling at regular intervals, or whether a certain amount of randomness in $\Delta_0$ is preferable. In the case of the MLE for the Ornstein–Uhlenbeck process, it can be seen from Table 1 that the randomness of the sampling scheme makes no difference for $\sigma^2$. On the other hand, for the Euler estimation of $\theta$ for the same process, randomness (i.e., higher $E[\Delta_0^2]$) adversely affects the bias but reduces the asymptotic variance. At the first order in $\varepsilon$, randomness has no effect on the estimation of $\sigma^2$; at the second order, more randomness reduces the asymptotic variance and the bias (since the first-order bias term is negative, a higher positive second-order bias works to reduce the bias).

Outside of the Ornstein–Uhlenbeck situation, it should be noted that even in the case of the MLE, it can occur that a somewhat random sampling can be preferable to sampling at a fixed interval. This occurs, for example, if one estimates $\sigma^2$ in the presence of a known drift function $\mu(x) = -x(1 - \exp(-x^4))$ (and, hence, known $\theta$). For that drift function, one then obtains that $E[(\partial^3\mu/\partial y^3)(Y_0)] > 0$ and so $\operatorname{sgn} \Omega_{\sigma^2}^{(3)} = -\operatorname{sgn} E[\Delta_0]$ since when we are only estimating $\sigma^2$, with $\theta_0$ known, the asymptotic variance of MLE is $\Omega_{\sigma^2} = \Omega^{(1)}_{\sigma^2}\varepsilon + \Omega^{(3)}_{\sigma^2}\varepsilon^3 + O(\varepsilon^4)$ with $\Omega^{(1)}_{\sigma^2} = 2\sigma_0^4 E[\Delta_0]$ and

$$(46) \qquad \Omega^{(3)}_{\sigma^2} = -\frac{1}{3}\sigma_0^6 E[\Delta_0]E[\Delta_0^2]E_{Y_0}\left[\frac{\partial^3 \mu}{\partial y^3}(Y_0;\theta_0)\right]$$

TABLE 1
*Asymptotic variance and bias for the Ornstein–Uhlenbeck process estimated using maximum likelihood and the Euler scheme. These expressions follow from specializing the general results to the Ornstein–Uhlenbeck process. When estimating $\theta$ with known $\sigma^2$ using the Euler scheme, $T_\theta = 0$ for the Ornstein–Uhlenbeck process because $h(Y_0, Y_1, \Delta, \bar{\beta}, \varepsilon)$ turns out to be a martingale. Note that it is perfectly acceptable for the variance of $\hat{\theta}$ to be below that of the MLE estimator. This can easily occur for an inconsistent estimator. Note that since $\theta_0 = \log(1 - \delta\bar{\theta})/\delta$, one can create a consistent estimator out of the Euler estimator $\hat{\theta}$ by using $\log(1 - \delta\hat{\theta})/\delta$. The latter is inefficient relative to the MLE estimator, as expected. When estimating $\sigma^2$ with known $\theta$, the first-order expansion for the MLE's $\Omega_{\sigma^2}$ is exact. This is because the Ornstein–Uhlenbeck process has a constant diffusion parameter and a Gaussian likelihood. But for $\theta$, the MLE's $\Omega_\theta$ involves an expansion because the exact log-likelihood of the process is a function of $\exp(-\theta\delta)$, which in our method is then Taylor-expanded in $\delta$*

|  | **MLE** | **Euler** |
|---|---|---|
| $\Omega_\theta$ | $2\theta_0 + \varepsilon^2\left(\frac{2\theta_0^3 E[\Delta_0^3]}{3E[\Delta_0]}\right) + O(\varepsilon^3)$ | $2\theta_0 - \varepsilon\left(\frac{2\theta_0^2 E[\Delta_0^2]}{E[\Delta_0]}\right) + O(\varepsilon^2)$ |
| $\bar{\theta} - \theta_0$ | $0$ | $-\varepsilon\left(\frac{\theta_0^2 E[\Delta_0^2]}{2E[\Delta_0]}\right) + O(\varepsilon^2)$ |
| $\Omega_{\sigma^2}$ | $\varepsilon(2\sigma_0^4 E[\Delta_0])$ | $\varepsilon(2\sigma_0^4 E[\Delta_0]) - \varepsilon^2(4\theta_0\sigma_0^4 E[\Delta_0]^2) + O(\varepsilon^3)$ |
| $\bar{\sigma}^2 - \sigma_0^2$ | $0$ | $-\varepsilon(\theta_0\sigma_0^2 E[\Delta_0]) + \varepsilon^2\left(\frac{2\theta_0^2 E[\Delta_0^2]}{3}\right) + O(\varepsilon^3)$ |

[see Aït-Sahalia and Mykland (2003) for an analysis of the MLE special case]. Since $\Omega^{(1)}_{\sigma^2}$ only depends on the first moment of $\Delta_0$, there is, therefore, a beneficial first-order effect of random sampling on $\Omega_{\sigma^2}$. For other drifts, such as, for instance, $\mu(x) = -x^3$, we have $E[(\partial^3\mu/\partial y^3)(Y_0)] < 0$ and, therefore, the opposite is true.

There is, therefore, no overall rule that covers all cases. In general, the impact of the sampling depends on the coefficients associated with the moments of $\Delta_0$, and the expansions derived in this paper can be used to gain insight into this impact.

## 5. Extensions of the theory.

5.1. *Extensions to more general estimating equations.* In terms of admissible $h$ functions, our theory can be extended from Taylor series to Laurent series (which have both positive and negative powers in $\varepsilon$). That is, the structure can be easily generalized to a situation where $h$ is of the form

$$h(y_1, y_0, \delta, \beta, \varepsilon) = \tilde{h}(y_1, y_0, \delta, \beta, \varepsilon) + \sum_{m=1}^{M} \frac{H_m(y_1, y_0, \delta, \beta, \varepsilon)}{\delta^m},$$

where $\tilde{h}$ and $\{H_m; m = 1, \ldots, M\}$ satisfy Assumption 3 with $\partial^k H_m(y_0, y_0, 0, \beta_0, 0)/\partial y_1^k = 0$ for $k = 1, \ldots, m$. Since this situation does not appear in practical estimation methods other than for $M = 1$, we have stated the result for that case, that is, (18), to avoid needlessly complicating the notation.

A different extension is the following. Instead of being of the form (18), the vector of moment functions $h$ is of the form

$$(47) \qquad h(y_1, y_0, \delta, \beta, \varepsilon) = \check{h}(y_1, y_0, \delta, \beta, \varepsilon) + \frac{K(y_1, y_0, \delta, \beta, \varepsilon)}{\varepsilon},$$

where both $\check{h}$ and $K$ can be Taylor expanded as specified by (55) and

$$(48) \quad K(y_0, y_0, 0, \beta_0, 0) = \frac{\partial K(y_0, y_0, 0, \beta_0, 0)}{\partial y_1} = \frac{\partial K(y_0, y_0, 0, \beta_0, 0)}{\partial \beta} = 0.$$

Then a simple modification of Lemmas 1 and 3 holds: evaluate (22) and (23) at $\check{h}$ instead of $h$, and replace (24), (25) and (33), respectively, by the following contributions from $K$:

$$(49) \qquad D^K_\beta = E_{\Delta, Y_0}[(\Gamma_{\beta_0} \cdot \dot{K})] + \frac{\varepsilon}{2}E_{\Delta, Y_0}[(\Gamma^2_{\beta_0} \cdot \dot{K})] + O(\varepsilon^2),$$

$$S_{\beta,0}^K = E_{\Delta,Y_0}[(\Gamma_{\beta_0} \cdot (\tilde{h} \times K'))] + \frac{\varepsilon}{2} E_{\Delta,Y_0}[(\Gamma_{\beta_0}^2 \cdot (\tilde{h} \times K'))]$$

$$+ E_{\Delta,Y_0}[(\Gamma_{\beta_0} \cdot (K \times \tilde{h}'))] + \frac{\varepsilon}{2} E_{\Delta,Y_0}[(\Gamma_{\beta_0}^2 \cdot (K \times \tilde{h}'))]$$

(50)
$$+ \frac{1}{2} E_{\Delta,Y_0}[(\Gamma_{\beta_0}^2 \cdot (K \times K'))]$$

$$+ \frac{\varepsilon}{6} E_{\Delta,Y_0}[(\Gamma_{\beta_0}^3 \cdot (K \times K'))] + O(\varepsilon^2),$$

$$[T_\beta^K]_{(i,j)} = \frac{1}{E[\Delta_0]} \bigg( \varepsilon^{\alpha_j - 1} E_{\Delta,Y_0}[(\Gamma_{\beta_0} \cdot K_i) \times r_j]$$

$$+ \frac{\varepsilon^{\alpha_j}}{2} E_{\Delta_0,Y_0}[(\Gamma_{\beta_0}^2 \cdot (K_i \times r_j))]$$

(51)
$$+ \varepsilon^{\alpha_i - 1} E_{\Delta_0,Y_0}[(\Gamma_{\beta_0} \cdot K_j) \times r_i]$$

$$+ \frac{\varepsilon^{\alpha_i}}{2} E_{\Delta_0,Y_0}[(\Gamma_{\beta_0}^2 \cdot (K_j \times r_i))] \bigg)$$

$$+ O(\varepsilon^{\min(\alpha_i,\alpha_j)+1}),$$

yielding $D_\beta = D_\beta^{\check{h}} + D_\beta^K$, $S_{\beta,0} = S_{\beta,0}^{\check{h}} + S_{\beta,0}^K$ and $T_\beta = T_\beta^{\check{h}} + T_\beta^K$.

Note that since $\varepsilon$ is deterministic, using $h$ or $\varepsilon h$ as the vector of moment functions produces the same estimator. Indeed, when $h$ is of the form (47), the two matrices $\Omega_\beta$ produced by applying Lemmas 1 and 3 with $(\tilde{h}, H) = (\varepsilon \check{h} + K, 0)$ or the first part of this remark with $(\check{h}, K)$ are identical.

5.2. *Extensions to more general Markov processes.* One can extend the theory to cover more general continuous-time Markov processes, such as jump-diffusions. In that case, the standard infinitesimal generator of the process applied to a smooth $f$ takes the form

$$J_{\beta_0} \cdot f = A_{\beta_0} \cdot f + \int \{f(y_1 + z, y_0, \delta, \beta, \varepsilon) - f(y_1, y_0, \delta, \beta, \varepsilon)\} \nu(dz, y_0),$$

where $A_{\beta_0}$, defined in (15), is the contribution coming from the diffusive part of the stochastic differential equation and $\nu(dz, y_0)$ is the Lévy jump measure specifying the number of jumps of size in $(z, z + dz)$ per unit of time [see, e.g., Protter (1992)]. In that case, our generalized infinitesimal generator becomes

$$\Gamma_{\beta_0} \cdot f \equiv \Delta_0 J_{\beta_0} \cdot f + \frac{\partial f}{\partial \varepsilon} + \frac{\partial f}{\partial \beta} \frac{\partial \beta}{\partial \varepsilon},$$

that is, the same expression as (17) except that $A_{\beta_0}$ is replaced by $J_{\beta_0}$.

5.3. *Extensions to more general sampling processes.* Another extension concerns the generation of the sampling intervals. For example, if the $\Delta_i$s are random and i.i.d., then $E[\Delta]$ has the usual meaning, but even if this is not the case, by $E[\Delta]$ we mean the limit (in probability, or just the limit if the $\Delta_i$s are nonrandom) of $\sum_{i=1}^{n} \Delta_i/n$ as $n$ tends to infinity. This permits the inclusion of the random non-i.i.d. and the nonrandom (but possibly irregularly spaced) cases for the $\Delta_i$s. At the cost of further complications, the theory can be extended to allow for dependence in the sampling intervals, whereby $\Delta_n$ is drawn conditionally on $(Y_{n-1}, \Delta_{n-1})$.

## 6. Proofs.

### 6.1. *Mixing.*

LEMMA 4. *Under Assumptions* 1 *and* 2, *the $\rho$-mixing coefficients of the discretely sampled process decay exponentially fast.*

6.2. *Proof of Lemma* 4. We start by showing that the sequence of $\rho$-mixing coefficients $\{\rho_\delta; \delta > 0\}$ of the process

$$(52) \qquad \rho_\delta \equiv \sup_{\{\phi,\psi \in L^2 | E[\phi(Y_0)] = E[\psi(Y_0)] = 0\}} \frac{E[\phi(Y_0)(U_\delta \cdot \psi)(Y_0)]}{\|\phi\|\|\psi\|},$$

decays exponentially fast as $\delta$ increases. Under Assumption 1, specifically condition (8), the operator $U_\delta$, as defined just after equation (28), is a strong contraction and there exists $\kappa > 0$ such that $\|U_\delta \cdot \psi\| \leq \exp(-\kappa\delta)\|\psi\|$ [see Propositions 8 and 9 in Hansen and Scheinkman (1995)]. Thus, by the Cauchy–Schwarz inequality,

$$|E[\phi(X_0)(U_\delta \cdot \psi)(X_0)]| \leq \|\phi\|\|U_\delta \cdot \psi\| \leq \|\phi\|\|\psi\|\exp(-\kappa\delta),$$

that is, $\rho_\delta \leq \exp(-\kappa\delta)$.

The mixing property of the underlying continuous time process $\{X_t; t \geq 0\}$ translates into the following mixing property for the discretely (and possibly randomly) sampled state process $\{Y_n; n = 0, \ldots, N_T\}$. For functions $\phi$ and $\psi$ in $L^2$, we have

$$\begin{aligned}
E[\phi(Y_0)\psi(Y_n)] &= E[\phi(X_0)\psi(X_{\Delta_1+\cdots+\Delta_n})] \\
&= E_{\Delta_1,\ldots,\Delta_n}[E[\phi(X_0)\psi(X_{\Delta_1+\cdots+\Delta_n})|\Delta_1,\ldots,\Delta_n]] \\
&= E_{\Delta_1,\ldots,\Delta_n}[E_{X_0}[\phi(X_0)E_{X_0}[\psi(X_{\Delta_1+\cdots+\Delta_n})|X_0,\Delta_1,\ldots,\Delta_n]]] \\
&= E_{\Delta_1,\ldots,\Delta_n}[E_{X_0}[\phi(X_0)(U_{\Delta_1+\cdots+\Delta_n} \cdot \psi)(X_0)]]
\end{aligned}$$

so that

$$(53) \qquad \begin{aligned}
|E[\phi(Y_0)\psi(Y_n)]| &\leq E_{\Delta_1,\ldots,\Delta_n}[|E_{X_0}[\phi(X_0)(U_{\Delta_1+\cdots+\Delta_n} \cdot \psi)(X_0)]|] \\
&\leq E_{\Delta_1,\ldots,\Delta_n}[\exp(-\lambda(\Delta_1+\cdots+\Delta_n))]\|\phi\|\|\psi\| \\
&= \{E_\Delta[\exp(-\kappa\Delta)]\}^n\|\phi\|\|\psi\|,
\end{aligned}$$

with the last equality following from the independence of the $\Delta_n$s. Since $0 < E_\Delta[\exp(-\kappa\Delta)] < 1$, the $Y_n$s satisfy a mixing property sufficient to insure the validity of the central limit theorem for sums of functions of the data $\{(\Delta_n, Y_n); n = 0, \ldots, N_T\}$.

6.3. *Proof of Lemma* 1. To calculate Taylor expansions of functions $f(Y_1, Y_0, \Delta, \bar{\beta}, \varepsilon) \in \mathcal{D}^J$, note first that

$$
\begin{aligned}
(54) \quad & E_{Y_1}[f(Y_1, Y_0, \Delta, \bar{\beta}, \varepsilon)|Y_0, \Delta = \delta] \\
&= \sum_{j=0}^{J} \frac{\varepsilon^j}{j!}(\Gamma^j_{\beta_0} \cdot f)(Y_0, Y_0, 0, \beta_0, 0) + O_p(\varepsilon^{J+1}).
\end{aligned}
$$

All the expectations are taken with respect to the law of the process at the true value $\beta_0$. This is in analogy to Theorem 1 in Aït-Sahalia and Mykland [(2003), page 498].

1. Starting with $D_\beta$, assume first that $H = 0$, and write a Taylor expansion of $E_{Y_1}[\dot{h}|Y_0, \Delta]$ in $\Delta$, using (54):

$$
\begin{aligned}
E_{\Delta, Y_1}&[\dot{h}(Y_1, Y_0, \Delta, \bar{\beta}, \varepsilon)|Y_0 = y_0] \\
&= \dot{h}(y_0, y_0, 0, \beta_0, 0) \\
&\quad + \varepsilon\left(E[\Delta_0][A_{\beta_0} \cdot \dot{h}] + \frac{\partial \dot{h}}{\partial \varepsilon} + \frac{\partial \dot{h}}{\partial \beta} \times \frac{\partial \bar{\beta}}{\partial \varepsilon}(\beta_0, 0)\right) + O(\varepsilon^2),
\end{aligned}
$$

with the partial derivatives on the right-hand side evaluated at $(y_0, y_0, 0, \beta_0, 0)$. This follows from the fact that $h$ can be Taylor expanded in $\varepsilon$ around 0,

$$
\begin{aligned}
(55) \quad h(y_1, y_0, \delta, \bar{\beta}, \varepsilon) &= h(y_0, y_0, 0, \beta_0, 0) + (y_1 - y_0)\frac{\partial h}{\partial y_1} + \frac{1}{2}(y_1 - y_0)^2 \frac{\partial^2 h}{\partial y_1^2} \\
&\quad + \frac{\partial h}{\partial \delta}\varepsilon\Delta_0 + \frac{\partial h}{\partial \varepsilon}\varepsilon + \frac{\partial h}{\partial \beta}\frac{\partial \bar{\beta}(\beta_0, 0)}{\partial \varepsilon}\varepsilon + o(\varepsilon),
\end{aligned}
$$

with all the partial derivatives of $h$ on the right-hand side evaluated at $(y_0, y_0, 0, \beta_0, 0)$. At the next order, we can write this more compactly as

$$
\begin{aligned}
(56) \quad E_{\Delta, Y_1}&[\dot{h}(Y_1, Y_0, \Delta, \bar{\beta}, \varepsilon)|Y_0 = y_0] \\
&= \dot{h}(y_0, y_0, 0, \beta_0, 0) + \varepsilon(\Gamma_{\beta_0} \cdot \dot{h})(Y_0, Y_0, 0, \beta_0, 0) \\
&\quad + \frac{\varepsilon^2}{2}(\Gamma^2_{\beta_0} \cdot \dot{h})(Y_0, Y_0, 0, \beta_0, 0) + O(\varepsilon^3).
\end{aligned}
$$

The unconditional expectation (22) follows from (56) by taking expectations with respect to $Y_0$ and using the law of iterated expectations.

Turning to $S_{\beta,0} \equiv E_{\Delta, Y_1, Y_0}[h(Y_1, Y_0, \Delta, \bar{\beta}, \varepsilon)h(Y_1, Y_0, \Delta, \bar{\beta}, \varepsilon)']$, assume first that $H = 0$. The result (23) follows from applying the generalized infinitesimal generator to $h \times h'$:

$E_{\Delta, Y_1}[(h \times h')(Y_1, Y_0, \Delta, \bar{\beta}, \varepsilon)|Y_0]$

$$= (h \times h')(Y_0, Y_0, 0, \beta_0, 0) + \varepsilon E_{\Delta, Y_0}[(\Gamma_{\beta_0} \cdot (h \times h'))(Y_0, Y_0, 0, \beta_0, 0)]$$
$$+ \frac{\varepsilon^2}{2} E_{\Delta, Y_0}[(\Gamma^2_{\beta_0} \cdot (h \times h'))(Y_0, Y_0, 0, \beta_0, 0)] + O_p(\varepsilon^3).$$

2. Suppose now that $H$ is not zero. Let $h_i = \tilde{h}_i + \Delta^{-1} H_i$ for $\tilde{h}_i \in \mathcal{D}^J$ and $H_i \in \mathcal{D}^{J+1}$. Applying (54) to $\tilde{h}_i$ and $H_i$ separately, then combining their expansions to get the expansion for $h_i$, we obtain that

$$\begin{aligned}
E_{Y_1}&[h_i(Y_1, Y_0, \Delta, \bar{\beta}, \varepsilon)|Y_0, \Delta] \\
&= E_{Y_1}[\tilde{h}_i(Y_1, Y_0, \Delta, \bar{\beta}, \varepsilon)|Y_0, \Delta] \\
&\quad + \varepsilon^{-1} E_{Y_1}[\Delta_0^{-1} H_i(Y_1, Y_0, \Delta, \bar{\beta}, \varepsilon)|Y_0, \Delta] \\
&= \sum_{j=0}^{J} \frac{\varepsilon^j}{j!} \Gamma^j_{\beta_0} \cdot \tilde{h}_i + \varepsilon^{-1} \sum_{j=1}^{J+1} \frac{\varepsilon^j}{j!} \Delta_0^{-1} \Gamma^j_{\beta_0} \cdot H_i + O_p(\varepsilon^{J+1}) \\
&= \sum_{j=0}^{J} \frac{\varepsilon^j}{j!} \left\{ (\Gamma^j_{\beta_0} \cdot \tilde{h}_i) + \frac{1}{j+1} \Delta_0^{-1}(\Gamma^{j+1}_{\beta_0} \cdot H_i) \right\} + O_p(\varepsilon^{J+1})
\end{aligned} \quad (57)$$

because under Assumption 3 we have $H_i = 0$ when evaluated at $(y_0, y_0, 0, \beta_0, 0)$. So the expansion (54) for $H_i$ starts at order $\varepsilon^1$ (or higher); without that, a singularity of order $\varepsilon^{-1}$ would result from the premultiplication by $\varepsilon^{-1}$.

The additional contribution to $D_\beta$ is given by (24) following a similar construction, where we use again equation (54). From

$$\begin{aligned}
E_{Y_1}&[\dot{H}(Y_1, Y_0, \Delta, \bar{\beta}, \varepsilon)|Y_0, \Delta] \\
&= \dot{H}(Y_0, Y_0, 0, \beta_0, 0) + \varepsilon(\Gamma_{\beta_0} \cdot \dot{H})(Y_0, Y_0, 0, \beta_0, 0) \\
&\quad + \frac{\varepsilon^2}{2}(\Gamma^2_{\beta_0} \cdot \dot{H})(Y_0, Y_0, 0, \beta_0, 0) + O_p(\varepsilon^3),
\end{aligned}$$

where we recall that $\dot{H}(Y_0, Y_0, 0, \beta_0, 0) = 0$ under (19) and

$$E_{\Delta, Y_1, Y_0}[\Delta^{-1} \dot{H}(Y_1, Y_0, \Delta, \bar{\beta}, \varepsilon)]$$
$$= E_{\Delta, Y_1, Y_0}[\Delta^{-1} E_{Y_1}[\dot{H}(Y_1, Y_0, \Delta, \bar{\beta}, \varepsilon)|Y_0, \Delta]],$$

we conclude that

$$\begin{aligned}
E_{\Delta, Y_1, Y_0}&[\Delta^{-1} \dot{H}(Y_1, Y_0, \Delta, \bar{\beta}, \varepsilon)] \\
&= E_{\Delta, Y_0}\left[ \Delta_0^{-1} \varepsilon^{-1} \left\{ \varepsilon(\Gamma_{\beta_0} \cdot \dot{H})(Y_0, Y_0, 0, \beta_0, 0) \right.\right. \\
&\qquad\qquad\qquad\qquad \left.\left. + \frac{\varepsilon^2}{2}(\Gamma^2_{\beta_0} \cdot \dot{H})(Y_0, Y_0, 0, \beta_0, 0) + O(\varepsilon^3) \right\} \right] \\
&= E_{\Delta, Y_0}[\Delta_0^{-1}(\Gamma_{\beta_0} \cdot \dot{H})(Y_0, Y_0, 0, \beta_0, 0)] \\
&\quad + \frac{\varepsilon}{2} E_{\Delta, Y_0}[\Delta_0^{-1}(\Gamma^2_{\beta_0} \cdot \dot{H})(Y_0, Y_0, 0, \beta_0, 0)] + O(\varepsilon^2).
\end{aligned}$$

The term contributed by $H$ to $S_{\beta,0}$ that is potentially the largest involves the cross product $(\Delta^{-1}H) \times (\Delta^{-1}H)$, that is, $E_{\Delta,Y_1,Y_0}[\Delta^{-2}(H \times H')(Y_1,Y_0,\Delta,\bar{\beta},\varepsilon)]$. To evaluate it, we start with

$$E_{Y_1}[(H \times H')(Y_1,Y_0,\Delta,\bar{\beta},\varepsilon)|Y_0,\Delta]$$
$$= (H \times H')(Y_0,Y_0,0,\beta_0,0) + \varepsilon(\Gamma_{\beta_0} \cdot (H \times H'))(Y_0,Y_0,0,\beta_0,0)$$
$$+ \frac{\varepsilon^2}{2}(\Gamma_{\beta_0}^2 \cdot (H \times H'))(Y_0,Y_0,0,\beta_0,0) + O_p(\varepsilon^3).$$

Next, note that

$$H(Y_0,Y_0,0,\beta_0,0) = 0 \quad \text{and} \quad (\Gamma_{\beta_0} \cdot (H \times H'))(Y_0,Y_0,0,\beta_0,0) = 0$$

under (19). Indeed, we have

$$(\Gamma_{\beta_0} \cdot (H \times H'))(Y_0,Y_0,0,\beta_0,0)$$
$$= \Delta_0 A_{\beta_0} \cdot (H \times H') + \frac{\partial(H \times H')}{\partial \varepsilon} + \frac{\partial(H \times H')}{\partial \bar{\beta}} \frac{\partial \bar{\beta}}{\partial \varepsilon}$$
$$= \Delta_0 \bigg\{ 2H \frac{\partial H'}{\partial \Delta} + \mu(Y_0;\theta_0) 2H \frac{\partial H'}{\partial y_1}$$
$$+ \sigma^2(Y_0;\gamma_0) \bigg( 2\frac{\partial H}{\partial y_1}\frac{\partial H'}{\partial y_1} + 2H\frac{\partial^2 H'}{\partial y_1 \partial y_1} \bigg) \bigg\} + 2H\frac{\partial H'}{\partial \varepsilon} + 2H\frac{\partial H'}{\partial \bar{\beta}}\frac{\partial \bar{\beta}}{\partial \varepsilon},$$

where in the equation above $H$ and its derivatives, listed without argument, are understood to be evaluated at $(Y_0,Y_0,0,\beta_0,0)$.

Since $H(Y_0,Y_0,0,\beta_0,0) = 0$ and $g_H = \partial H(Y_0,Y_0,0,\beta_0,0)/\partial y_1 = 0$ under (19), it follows that

$$(\Gamma_{\beta_0} \cdot (H \times H'))(Y_0,Y_0,0,\beta_0,0) = 0.$$

Then, from

$$E_{\Delta,Y_1,Y_0}[\Delta^{-2}(H \times H')(Y_1,Y_0,\Delta,\bar{\beta},\varepsilon)]$$
$$= E_{\Delta,Y_1,Y_0}[\Delta^{-2}E_{Y_1}[(H \times H')(Y_1,Y_0,\Delta,\bar{\beta},\varepsilon)|Y_0,\Delta]],$$

we conclude that

$$E_{\Delta,Y_1,Y_0}[\Delta^{-2}(H \times H')(Y_1,Y_0,\Delta,\bar{\beta},\varepsilon)]$$
$$= E_{\Delta,Y_0}\bigg[\Delta_0^{-2}\varepsilon^{-2}\bigg\{\frac{\varepsilon^2}{2}(\Gamma_{\beta_0}^2 \cdot (H \times H'))(Y_0,Y_0,0,\beta_0,0)$$
$$+ \frac{\varepsilon^3}{6}(\Gamma_{\beta_0}^3 \cdot (H \times H'))(Y_0,Y_0,0,\beta_0,0) + O(\varepsilon^4)\bigg\}\bigg]$$
$$= E_{\Delta,Y_0}[\Delta_0^{-2}(\Gamma_{\beta_0}^2 \cdot (H \times H'))(Y_0,Y_0,0,\beta_0,0)]$$
$$+ \frac{\varepsilon}{6}E_{\Delta,Y_0}[\Delta_0^{-2}(\Gamma_{\beta_0}^3 \cdot (H \times H'))(Y_0,Y_0,0,\beta_0,0)] + O(\varepsilon^2).$$

Finally, the other two cross product terms, $E_{\Delta,Y_1,Y_0}[\Delta^{-1}(H \times \tilde{h}')]$ and $E_{\Delta,Y_1,Y_0}[\Delta^{-1}(\tilde{h} \times H')]$, are dealt with similarly. They are of order $O(1)$ since $H(Y_0,Y_0,0,\beta_0,0) = 0$:

$$\begin{aligned}
E_{Y_1}[\Delta^{-1}&(H \times \tilde{h}')(Y_1,Y_0,\Delta,\bar{\beta},\varepsilon)|Y_0,\Delta] \\
&= \varepsilon^{-1}\Delta_0^{-1} E_{Y_1}[(H \times \tilde{h}')(Y_1,Y_0,\Delta,\bar{\beta},\varepsilon)|Y_0,\Delta] \\
&= \varepsilon^{-1}\Delta_0^{-1}\left\{(H \times \tilde{h}') + \varepsilon(\Gamma_{\beta_0} \cdot (H \times \tilde{h}')) + \frac{\varepsilon^2}{2}(\Gamma_{\beta_0}^2 \cdot (H \times \tilde{h}')) + O_p(\varepsilon^3)\right\} \\
&= \Delta_0^{-1}(\Gamma_{\beta_0} \cdot (H \times \tilde{h}')) + \frac{\varepsilon}{2}\Delta_0^{-1}(\Gamma_{\beta_0}^2 \cdot (H \times \tilde{h}')) + O_p(\varepsilon^3)
\end{aligned}$$

and similarly for $E_{Y_1}[\Delta^{-1}(\tilde{h} \times H')(Y_1,Y_0,\Delta,\bar{\beta},\varepsilon)|Y_0,\Delta]$.

6.4. *Proof of Lemma* 2. Note first that $r_i(y,\beta_0,\varepsilon)$ and $\check{r}_i(y,\beta_0,\varepsilon)$ are well defined as a consequence of the $L^2$ boundedness of $A_{\beta_0} \cdot q_i$, and the exponential mixing from that follows from Lemma 4. We here take the first expression in (30) to be the definition of $\check{r}$. To see the equality with the second expression, note that $A_{\beta_0} \cdot (tq_i) = q_i + tA_{\beta_0} \cdot q_i$. The second expression for $\check{r}$ follows. As before, $Y$ has the stationary distribution of $X_0$.

Let $N_0(u)$ be the number of $\tau_i^{(0)} = \tau_i/\varepsilon$ in the interval $(0,u]$. Also, set $Z(t) = E[\tau_{N(t)+1} - t]$ and $Z^{(0)}(t) = E[\tau_{N_0(t)+1}^{(0)} - t]$, and note that, by Wald's identity,

$$(58) \qquad Z^{(0)}(u) = E[\Delta_0]E[N^{(0)}(u) + 1] - u,$$

and similarly without the superscript 0. In particular, $Z(t) = \varepsilon Z^{(0)}(t/\varepsilon)$. Since the integrals are well defined, it follows that

$$(59) \quad \begin{aligned}\varepsilon^{-1}(r_i(Y,\beta_0,\varepsilon) - \check{r}_i(Y,\beta_0,\varepsilon)) &= \varepsilon^{-1}\int_0^\infty U_t \cdot A_{\beta_0} \cdot q_i(Y,\beta_0,\varepsilon)Z(t)\,dt \\ &= -\int_0^\infty U_t \cdot A_{\beta_0} \cdot q_i(Y,\beta_0,\varepsilon)Z^{(0)}(t/\varepsilon)\,dt.\end{aligned}$$

In the sequel, we assume that $\varepsilon \to 0$ through a countable sequence. The $L^2$ limit will be independent of the choice of sequence, and so it will be valid as $\varepsilon$ goes to zero generally. We also need the mixing coefficient $\lambda$ from Lemma 4 (there written as $\kappa$) and an exponent $\lambda_1 > 0$ which can take on different values.

We first need to establish some facts about $Z^{(0)}(t)$, and here we make use of Feller (1971), to which all references in the next two paragraphs are made. First note that $Z^{(0)}$ is the solution of the renewal equation $Z^{(0)} = z^{(0)} + F^{(0)} * Z^{(0)}$, where $F^{(0)}$ is the c.d.f. of $\Delta_0$, and $z^{(0)}(t) = \int_t^\infty (1 - F^{(0)}(\delta))\,d\delta$. This follows from the proof of Theorem XI-3.1 (pages 366 and 367). Since we have assumed that $E[\Delta_0^2] < \infty$, the same proof assures that $\limsup_t Z^{(0)}(t) < \infty$ in the nonarithmetic case for $\Delta_0$, and the same follows in the arithmetic case from the development on pages 362 and 363. (The distinction between the arithmetic and nonarithmetic cases is described on page 138.) Since $z^{(0)}$ is bounded, the Lemma on page 359 assures that $Z^{(0)}$ is bounded on finite intervals, whence

$$\text{(60)} \qquad \sup_t Z^{(0)}(t) < \infty \quad \text{and} \quad \inf_t Z^{(0)}(t) \geq 0,$$

where the latter inequality is by construction.

Also, the same Theorem XI-3.1 in Feller (1971) establishes that $Z^{(0)}(t) \to \frac{1}{2}E[\Delta_0^2]/E[\Delta_0]$ as $t \to \infty$ in the nonarithmetic case. In this case, therefore, for all $\lambda_1 > 0$, in the sense of weak convergence of measures on $[0, \infty)$,

$$\text{(61)} \qquad \exp\{-\lambda_1 t\} Z^{(0)}\left(\frac{t}{\varepsilon}\right) dt \to \exp\{-\lambda_1 t\} \frac{1}{2} \frac{E[\Delta_0^2]}{E[\Delta_0]} dt,$$

by (60). In the arithmetic case, $Z^{(0)}(t)$ does not converge, but (61) follows from the results on pages 362 and 363. This is what we needed from Feller (1971), and we now proceed to make use of (60) and (61).

We then establish the convergence in probability of (59). As in the proof of Lemma 4,

$$\|U_t(A_{\beta_0} \cdot q_i(Y, \beta_0, \varepsilon) - A_{\beta_0} \cdot q_i(Y, \beta_0, 0))\|$$
$$\leq \exp\{-\lambda t\} \|A_{\beta_0} \cdot q_i(Y, \beta_0, \varepsilon) - A_{\beta_0} \cdot q_i(Y, \beta_0, 0)\|.$$

By the $L^2$ continuity of $A_{\beta_0} \cdot q_i(Y, \beta_0, \varepsilon)$, and by (60), we can replace $U_t A_{\beta_0} \cdot q_i(Y, \beta_0, \varepsilon)$ by $U_t A_{\beta_0} \cdot q_i(Y, \beta_0, 0)$ for the purpose of this convergence. Since $U_t A_{\beta_0} \cdot q_i(Y, \beta_0, 0)$ can be taken to be continuous in $t$ on $[0, \infty]$ (since the limit is zero as $t \to \infty$), and in view of (61) (with $\lambda_1 < \lambda$), the limit of (59) must be as in (29), but for the moment we have only shown convergence in probability.

The final result (29) and (30) then follows if we can show that the square of the left-hand side of (59) is uniformly integrable as $\varepsilon \to 0$. This is the case since

$$E[\varepsilon^{-2}(r_i(Y, \beta_0, \varepsilon) - \breve{r}_i(Y, \beta_0, \varepsilon))^2]$$
$$= \int_0^\infty dt \int_0^\infty ds E[U_t \cdot A_{\beta_0} \cdot q_i(Y, \beta_0, \varepsilon) U_s \cdot A_{\beta_0} \cdot q_i(Y, \beta_0, \varepsilon)]$$
$$\times Z^{(0)}(t/\varepsilon) Z^{(0)}(s/\varepsilon).$$

In the same way as in the discussion above, the limit of the integral coincides with the integral of the limit. Hence, uniform integrability follows.

To see how $\breve{r}$ solves the differential equation, with the given side condition, proceed as follows. By the second expression in (30), and since $A_{\beta_0}$ and $U_t$ commute,

$$(A_{\beta_0} \cdot \breve{r}_i)(y, \beta_0, \varepsilon) = \int_0^\infty (U_t \cdot A_{\beta_0} \cdot q_i)(y, \beta_0, \varepsilon) = -q_i(y, \beta_0, \varepsilon).$$

If $\breve{r}_i$ is chosen to satisfy

(62) $$E_{Y_0}[\breve{r}_i(Y_0, \beta_0, \varepsilon)] = 0$$

under the stationary distribution, asymptotic ergodicity will force $\breve{r}_i$ to have the second form from (30).

Exploiting the form of the scale function $s$ defined in (5), we can rewrite (31) as

$$\frac{\partial}{\partial y}\left[\frac{\partial \breve{r}_i(y, \beta_0, \varepsilon)}{\partial y} \frac{1}{s(y; \beta_0)}\right]$$
$$= \frac{\partial}{\partial y}\left[\frac{1}{s(y; \beta_0)}\right]\frac{\partial \breve{r}_i(y, \beta_0, \varepsilon)}{\partial y} + \frac{\partial^2 \breve{r}_i(y, \beta_0, \varepsilon)}{\partial y^2}\frac{1}{s(y; \beta_0)}$$
$$= \left(2\frac{\mu(y; \theta_0)}{\sigma^2(y; \gamma_0)}\frac{\partial \breve{r}_i(y, \beta_0, \varepsilon)}{\partial y} + \frac{\partial^2 \breve{r}_i(y, \beta_0, \varepsilon)}{\partial y^2}\right)\frac{1}{s(y; \beta_0)}$$
$$= -\frac{2\, q_i(y, \beta_0, \varepsilon)}{\sigma^2(y; \gamma_0) s(y; \beta_0)}.$$

To solve this, we have

(63) $$\frac{\partial \breve{r}_i(y, \beta_0, \varepsilon)}{\partial y} = s(y; \beta_0)\left(C_1 - \int_{\underline{x}}^y \frac{2\, q_i(x, \beta_0, \varepsilon)}{\sigma^2(x; \gamma_0) s(x; \beta_0)}\, dx\right).$$

Subject to regularity conditions on the function $\sigma^2$, the constant of integration must be $C_1 = 0$, otherwise $\breve{r}_i$ would not be integrable under $\pi$. It follows that

(64) $$\breve{r}_i(y, \beta_0, \varepsilon) = C_2 - \int_{\underline{x}}^y \int_{\underline{x}}^z \frac{2\, q_i(x, \beta_0, \varepsilon)}{\sigma^2(x; \gamma_0) s(x; \beta_0)}\, dx\, s(z, \beta_0)\, dz,$$

where the second constant of integration $C_2$ is determined so that (62) holds. We only need the function $\breve{r}$ for the purpose of calculating expressions of the form $E_{Y_0}[\phi(Y_0)\breve{r}_i(Y_0, \beta_0, \varepsilon)]$, where $E_{Y_0}[\phi(Y_0)] = 0$ (as when $\phi = q$, for instance). Then the value of $C_2$ is irrelevant for the calculation of those unconditional expectations.

As $\varepsilon \to 0$, we have $r_i(y, \beta_0, 0) = \breve{r}_i(y, \beta_0, 0)$ and it follows from (63) that

(65) $$\frac{\partial}{\partial y}\left[\frac{\partial r_i(y, \beta_0, 0)}{\partial y} \frac{1}{s(y; \beta_0)}\right] = -\frac{2\, q_i(y, \beta_0, 0)}{\sigma^2(y; \gamma_0) s(y; \beta_0)}$$

since that equation does not involve differentiation with respect to $\varepsilon$. Indeed, in light of (29), we define $\partial r_i/\partial \varepsilon$ as follows:

$$\frac{\partial r_i}{\partial \varepsilon}(y, \beta_0, 0) = \frac{\partial \tilde{r}_i}{\partial \varepsilon}(y, \beta_0, 0) + \frac{1}{2}\frac{E[\Delta_0^2]}{E[\Delta_0]}q_i(y, \beta_0, 0).$$

We also define
$$\frac{\partial^k r_i(Y_0, \beta_0, 0)}{\partial y^k} \equiv \frac{\partial^k \check{r}_i(Y_0, \beta_0, 0)}{\partial y^k}$$
for $k = 1, 2$, and with these definitions of the partial derivatives of $r_i$ evaluated at $(Y_0, \beta_0, 0)$ we see that $r_i$ is Taylor-expandable in the form
$$r_i(Y_1, \beta_0, \varepsilon) = r_i(Y_0, \beta_0, 0) + (Y_1 - Y_0)\frac{\partial r_i(Y_0, \beta_0, 0)}{\partial y}$$
$$+ \frac{1}{2}(Y_1 - Y_0)^2 \frac{\partial^2 r_i(Y_0, \beta_0, 0)}{\partial y^2} + \varepsilon \frac{\partial r_i(Y_0, \beta_0, 0)}{\partial \varepsilon} + o_p(\varepsilon).$$

If $\sigma^2 = \gamma$ constant, dividing (65) by $\sigma_0^2$ yields an equivalent form in terms of the stationary density $\pi$:
$$\frac{\partial}{\partial y}\left[\frac{\partial r_i(y, \beta_0, 0)}{\partial y}\pi(y; \beta_0)\right] = -\frac{2}{\sigma_0^2}q_i(y, \beta_0, 0)\pi(y; \beta_0).$$

6.5. *Proof of Lemma 3.* 1. When the moment condition is not a martingale, the matrix $S_\beta$ includes time series terms $T_\beta = S_\beta - S_{\beta,0}$ which must be calculated. We start by showing the derivation in the case of scalar $h$; the generalization to the vector case is straightforward and is given at the end of this part of the proof. Recall equation (27), now for a scalar, $E_{\Delta, Y_1}[h(Y_1, Y_0, \Delta, \bar{\beta}, \varepsilon)|Y_0] = \varepsilon^\alpha q(Y_0, \beta_0, \varepsilon)$, where $q(Y_0, \beta_0, \varepsilon)$ is of order $O(1)$ in $\varepsilon$, and where the $\alpha$ is an integer greater than zero, typically $\alpha = 1$ or $2$. The covariance terms then become

$$
\begin{aligned}
T_\beta = S_\beta - S_{\beta,0} &= 2\sum_{k=1}^\infty S_{\beta,k} \\
&= 2\sum_{j=1}^\infty E[h(Y_1, Y_0, \Delta^{(0)}, \bar{\beta}, \varepsilon)h(Y_{k+1}, Y_k, \Delta^{(k)}, \bar{\beta}, \varepsilon)] \\
&= 2\sum_{k=1}^\infty E[h(Y_1, Y_0, \Delta^{(0)}, \bar{\beta}, \varepsilon)E[h(Y_{k+1}, Y_k, \Delta^{(k)}, \bar{\beta}, \varepsilon)|\Im_k]] \\
&= 2\sum_{k=1}^\infty E[h(Y_1, Y_0, \Delta^{(0)}, \bar{\beta}, \varepsilon)\varepsilon^\alpha q(Y_k, \beta_0, \varepsilon)] \\
&= 2\varepsilon^\alpha \sum_{k=1}^\infty E[h(Y_1, Y_0, \Delta^{(0)}, \bar{\beta}, \varepsilon)E[q(Y_k, \beta_0, \varepsilon)|Y_1]] \\
&= 2\varepsilon^{\alpha-1}\frac{1}{E[\Delta_0]}E_{\Delta, Y_1, Y_0}[h(Y_1, Y_0, \Delta, \bar{\beta}, \varepsilon)r(Y_1, \beta_0, \varepsilon)],
\end{aligned}
$$
(66)

where $\Im_j$ denotes the standard filtration up to time $j$.

The final transition in (66) requires showing that

(67) $$r(y, \beta_0, \varepsilon) = \varepsilon E[\Delta_0]\sum_{k=1}^\infty E_{Y_k}[q(Y_k, \beta_0, \varepsilon)|Y_1 = y].$$

To see this, note that $q(X_0, \beta_0, \varepsilon)$ and $(A_{\beta_0} \cdot q)(X_0, \beta_0, \varepsilon)$ are integrable under the stationary distribution, for $\varepsilon \in [0, \varepsilon_0]$. Then, for $t \geq u$,

$$E\left[\int_{\tau_{n-1} \wedge u}^{t} (A_{\beta_0} \cdot q)(X_s, \varepsilon) \, ds \Big| X_0 = y\right]$$
$$= E\left[\int_{0}^{t} (A_{\beta_0} \cdot q)(X_s, \varepsilon) I_{(s \geq \tau_{n-1} \wedge u)} \, ds \Big| X_0 = y\right]$$
$$= \int_{0}^{t} E[(A_{\beta_0} \cdot q)(X_s, \varepsilon) | X_0 = y] P(s \geq \tau_{n-1} \wedge u) \, ds$$
$$= \int_{0}^{t} (U_s A_{\beta_0} \cdot q)(y, \varepsilon) P(s \geq \tau_{n-1} \wedge u) \, ds.$$

The validity of Fubini's theorem and the integrability of all quantities considered follow from our assumptions since also $\tau_{n-1}$ is independent of the $X$ process, and the latter is stationary. These facts are also used in the following.

By Itô's lemma, and since $\int_{\tau_{n-1} \wedge u}^{t} \frac{\partial}{\partial y} q(X_s, \varepsilon) \sigma(X_s; \gamma_0) \, dW_s$ is a local martingale in $t$, we therefore get

$$(68) \quad E[q(Y_{n-1}, \varepsilon) | X_0 = y] = -\int_{0}^{+\infty} (U_s A_{\beta_0} \cdot q)(y, \varepsilon) P(s \geq \tau_{n-1}) \, ds.$$

This is by first letting $t \to +\infty$ and then $u \to +\infty$. We here use that $E[q(X_t, \varepsilon)]$ goes to zero as $t$ gets large.

To go from (68) to (67), note that the former implies

$$(69) \quad \begin{aligned} &\varepsilon E[\Delta_0] \sum_{k=1}^{n} E_{Y_k}[q(Y_k, \beta_0, \varepsilon) | Y_1 = y] \\ &= -\varepsilon E[\Delta_0] \int_{0}^{+\infty} (U_s A_{\beta_0} \cdot q)(y, \varepsilon) \left( \sum_{k=1}^{n} P(s \geq \tau_{k-1}) \right) ds. \end{aligned}$$

As $n \to +\infty$, we have $\sum_{k=1}^{n} P(s \geq \tau_{k-1}) \to E[N_s] + 1$. Note that $E[N_s] < +\infty$ by the Lemma on page II-359 in Feller (1971). Also, since $E[\Delta] < +\infty$, $E[\tau_{N_s+1}] = E[\Delta](E[N_s] + 1)$. It follows that one can let $n$ go to infinity in (69) and still have a finite limit. The result (67) follows.

We now proceed with the analysis of $T_\beta$. Assume first that $H = 0$. We return to the general case below. From (66),

$$T_\beta = 2\varepsilon^{\alpha-1} \frac{1}{E[\Delta_0]} E_{\Delta, Y_1, Y_0}[h(Y_1, Y_0, \Delta, \bar{\beta}, \varepsilon) r(Y_1, \beta_0, \varepsilon)]$$
$$= 2\varepsilon^{\alpha-1} \frac{1}{E[\Delta_0]} (E_{Y_0}[h(Y_0, Y_0, 0, \beta_0, 0) r(Y_0, \beta_0, 0)]$$
$$+ \varepsilon E_{\Delta_0, Y_0}[(\Gamma_{\beta_0} \cdot (h \times r))(Y_0, Y_0, 0, \beta_0, 0)]$$

$$+ \frac{\varepsilon^2}{2} E_{\Delta_0, Y_0}[(\Gamma_{\beta_0}^2 \cdot (h \times r))(Y_0, Y_0, 0, \beta_0, 0)] + O_p(\varepsilon^3))$$

$$= \frac{2}{E[\Delta_0]} (\varepsilon^{\alpha-1} E_{Y_0}[(h \times r)(Y_0, Y_0, 0, \beta_0, 0)]$$

$$+ \varepsilon^\alpha E_{\Delta_0, Y_0}[(\Gamma_{\beta_0} \cdot (h \times r))(Y_0, Y_0, 0, \beta_0, 0)]$$

$$+ \frac{\varepsilon^{\alpha+1}}{2} E_{\Delta_0, Y_0}[(\Gamma_{\beta_0}^2 \cdot (h \times r))(Y_0, Y_0, 0, \beta_0, 0)])$$

$$+ O(\varepsilon^{\alpha+2}),$$

where $(h \times r)(Y_0, Y_0, 0, \beta_0, 0) \equiv h(Y_0, Y_0, 0, \beta_0, 0) \, r(Y_0, \beta_0, 0)$, and

$$(\Gamma_{\beta_0} \cdot (h \times r))(Y_0, Y_0, 0, \beta_0, 0)$$

$$= (\Gamma_{\beta_0} \cdot h) \times r + h \times (\Gamma_{\beta_0} \cdot r) + \Delta_0 \sigma^2(Y_0; \gamma_0) \frac{\partial r}{\partial y_1} \frac{\partial h}{\partial y_1}$$

$$(70) \quad = \left(\Delta_0 \left(\frac{\partial h}{\partial \Delta} + \mu(Y_0; \theta_0) \frac{\partial h}{\partial y_1} + \frac{\sigma^2(Y_0; \gamma_0)}{2} \frac{\partial^2 h}{\partial y_1^2}\right) + \frac{\partial h}{\partial \varepsilon} + \frac{\partial h}{\partial \beta} \frac{\partial \beta}{\partial \varepsilon}\right) \times r$$

$$+ h \times \left(\Delta_0 \left(\mu(Y_0; \theta_0) \frac{\partial r}{\partial y_1} + \frac{\sigma^2(Y_0; \gamma_0)}{2} \frac{\partial^2 r}{\partial y_1^2}\right) + \frac{\partial r}{\partial \varepsilon}\right)$$

$$+ \Delta_0 \sigma^2(Y_0; \gamma_0) \frac{\partial h}{\partial y_1} \times \frac{\partial r}{\partial y_1},$$

with the understanding here and below that the functions listed without arguments are all evaluated at $Y_1 = Y_0$, $\Delta = 0$, $\bar{\beta} = \beta_0$ [since $\bar{\beta}(\beta_0, 0) = \beta_0$] and $\varepsilon = 0$.

Note that this requires that the function $r$ be Taylor-expandable in $\varepsilon$ as given in (66).

For multidimensional $h = (h_1, \ldots, h_r)'$, still assuming $H = 0$, the $(i,j)$ term of the $T_\beta = S_\beta - S_{\beta,0}$ matrix is

$$(71) \quad \begin{aligned}[T_\beta]_{(i,j)} &= \sum_{k=1}^\infty \{E[h_i(Y_1, Y_0, \Delta^{(0)}, \bar{\beta}, \varepsilon) h_j(Y_{k+1}, Y_k, \Delta^{(k)}, \bar{\beta}, \varepsilon)] \\ &\quad + E[h_j(Y_1, Y_0, \Delta^{(0)}, \bar{\beta}, \varepsilon) h_i(Y_{k+1}, Y_k, \Delta^{(k)}, \bar{\beta}, \varepsilon)]\} \\ &= \varepsilon^{\alpha_j - 1} \frac{1}{E[\Delta_0]} E_{\Delta, Y_1, Y_0}[h_i(Y_1, Y_0, \Delta, \bar{\beta}, \varepsilon) r_j(Y_1, \beta_0, \varepsilon)] \\ &\quad + \varepsilon^{\alpha_i - 1} \frac{1}{E[\Delta_0]} E_{\Delta, Y_1, Y_0}[h_j(Y_1, Y_0, \Delta, \bar{\beta}, \varepsilon) r_i(Y_1, \beta_0, \varepsilon)].\end{aligned}$$

By applying the univariate calculation above to the two terms involving $h_i$ and $h_j$, it follows that $[T_\beta]_{(i,j)}$ is given by

$$[T_\beta]_{(i,j)} = \frac{1}{E[\Delta_0]} (\varepsilon^{\alpha_j - 1} E_{Y_0}[(h_i \times r_j)(Y_0, Y_0, 0, \beta_0, 0)]$$

$$+ \varepsilon^{\alpha_j} E_{\Delta_0, Y_0}[(\Gamma_{\beta_0} \cdot (h_i \times r_j))(Y_0, Y_0, 0, \beta_0, 0)]$$

$$+ \frac{\varepsilon^{\alpha_j+1}}{2} E_{\Delta_0, Y_0}[(\Gamma_{\beta_0}^2 \cdot (h_i \times r_j))(Y_0, Y_0, 0, \beta_0, 0)])$$

$$+ \frac{1}{E[\Delta_0]} (\varepsilon^{\alpha_i-1} E_{Y_0}[(h_j \times r_i)(Y_0, Y_0, 0, \beta_0, 0)]$$

$$+ \varepsilon^{\alpha_i} E_{\Delta_0, Y_0}[(\Gamma_{\beta_0} \cdot (h_j \times r_i))(Y_0, Y_0, 0, \beta_0, 0)]$$

$$+ \frac{\varepsilon^{\alpha_i+1}}{2} E_{\Delta_0, Y_0}[(\Gamma_{\beta_0}^2 \cdot (h_j \times r_i))(Y_0, Y_0, 0, \beta_0, 0)])$$

$$+ O(\varepsilon^{\min(\alpha_i, \alpha_j)+2}).$$

2. We now investigate the contribution of a nonzero $H$ to $T_\beta$. Equation (27) now follows from

$$\begin{aligned}
& E_{\Delta, Y_1}[h_i(Y_1, Y_0, \Delta, \bar{\beta}, \varepsilon)|Y_0] \\
&= E_\Delta[E_{Y_1}[h_i(Y_1, Y_0, \Delta, \bar{\beta}, \varepsilon)|Y_0, \Delta]] \\
&= \sum_{j=0}^{J} \frac{\varepsilon^j}{j!} \left\{ E_{\Delta_0}[(\Gamma_{\beta_0}^j \cdot \tilde{h}_i)] + \frac{1}{(j+1)} E_{\Delta_0}[\Delta_0^{-1}(\Gamma_{\beta_0}^{j+1} \cdot H_i)] \right\} + O(\varepsilon^{J+1}) \\
&= \varepsilon^{\alpha_i} q_i(Y_0, \beta_0, 0) + O_p(\varepsilon^{\alpha_i+1})
\end{aligned} \tag{72}$$

if we let $\alpha_i$ denote an index $j$ at which the sum in the right-hand side of (72) is nonzero. As above, consider first the case of scalar $H$ and recall that $h = \tilde{h} + \Delta^{-1} H$. We now have to look at

$$\begin{aligned}
T_\beta &= 2 \sum_{j=1}^{\infty} E[h(Y_1, Y_0, \Delta^{(0)}, \bar{\beta}, \varepsilon) h(Y_{k+1}, Y_k, \Delta^{(k)}, \bar{\beta}, \varepsilon)] \\
&= 2 \sum_{k=1}^{\infty} E[\{\tilde{h}(Y_1, Y_0, \Delta^{(0)}, \bar{\beta}, \varepsilon) \\
&\qquad + (\Delta^{(0)})^{-1} H(Y_1, Y_0, \Delta^{(0)}, \bar{\beta}, \varepsilon)\} E[h(Y_{k+1}, Y_k, \Delta^{(k)}, \bar{\beta}, \varepsilon)|\Im_k]] \\
&= 2\varepsilon^{\alpha-1} \frac{1}{E[\Delta_0]} E_{\Delta, Y_1, Y_0}[\{\tilde{h}(Y_1, Y_0, \Delta, \bar{\beta}, \varepsilon) \\
&\qquad + \Delta^{-1} H(Y_1, Y_0, \Delta, \bar{\beta}, \varepsilon)\} r(Y_1, \beta_0, \varepsilon)] \\
&\quad + O(\varepsilon^\alpha),
\end{aligned}$$

where the term $E_{\Delta, Y_1, Y_0}[\tilde{h}(Y_1, Y_0, \Delta, \bar{\beta}, \varepsilon) r(Y_1, \beta_0, \varepsilon)]$ is the one we dealt with in part 1 of this proof. The additional contribution to $T_\beta$ is, therefore, rep-

resented by the term

$$
\begin{aligned}
(73)\quad T_\beta^H &= 2\varepsilon^{\alpha-1}\frac{1}{E[\Delta_0]}E_{\Delta,Y_1,Y_0}[\Delta^{-1}H(Y_1,Y_0,\Delta,\bar\beta,\varepsilon)r(Y_1,\beta_0,\varepsilon)] \\
&= 2\varepsilon^{\alpha-1}\frac{1}{E[\Delta_0]}E_{\Delta,Y_0}[\Delta^{-1}E_{Y_1}[H(Y_1,Y_0,\Delta,\bar\beta,\varepsilon)r(Y_1,\beta_0,\varepsilon)|Y_0,\Delta]].
\end{aligned}
$$

By (54), the conditional expectation of $H \times r$ can be Taylor-expanded as

$$
\begin{aligned}
(74)\quad &E_{Y_1}[H(Y_1,Y_0,\Delta,\bar\beta,\varepsilon)r(Y_1,\beta_0,\varepsilon)|Y_0,\Delta] \\
&= H(Y_0,Y_0,0,\beta_0,0)r(Y_0,\beta_0,0) \\
&\quad + \varepsilon(\Gamma_{\beta_0}\cdot(H\times r))(Y_0,Y_0,0,\beta_0,0) \\
&\quad + \frac{\varepsilon^2}{2}(\Gamma^2_{\beta_0}\cdot(H\times r))(Y_0,Y_0,0,\beta_0,0) + O_p(\varepsilon^3).
\end{aligned}
$$

Recall that under (19), $H(Y_0,Y_0,0,\beta_0,0)=0$ so the term of order $\varepsilon^0$ in (74) is 0. For the term of order $\varepsilon^1$, we have as in (70),

$$
\begin{aligned}
(\Gamma_{\beta_0}\cdot(H\times r))&(Y_0,Y_0,0,\beta_0,0) \\
&= (\Gamma_{\beta_0}\cdot H)\times r + H\times(\Gamma_{\beta_0}\cdot r) + \Delta_0\sigma^2(Y_0;\gamma_0)\frac{\partial r}{\partial y_1}\frac{\partial H}{\partial y_1} \\
&= (\Gamma_{\beta_0}\cdot H)\times r,
\end{aligned}
$$

with the last equation following from the fact that

$$H(Y_0,Y_0,0,\beta_0,0) = \frac{\partial H}{\partial y_1}(Y_0,Y_0,0,\beta_0,0) = 0$$

under (19). Next,

$$
\begin{aligned}
E_{\Delta,Y_0}&[\Delta^{-1}\varepsilon(\Gamma_{\beta_0}\cdot(H\times r))(Y_0,Y_0,0,\beta_0,0)] \\
&= E_{\Delta,Y_0}[\Delta_0^{-1}(\Gamma_{\beta_0}\cdot H)\times r] \\
&= E_{\Delta,Y_0}\left[\Delta_0^{-1}\left(\Delta_0\left(\frac{\partial H}{\partial\Delta}+\mu(Y_0;\theta_0)\frac{\partial H}{\partial y_1}+\frac{\sigma^2(Y_0;\gamma_0)}{2}\frac{\partial^2 H}{\partial y_1^2}\right)\right.\right. \\
&\qquad\qquad\qquad\left.\left. + \frac{\partial H}{\partial\varepsilon}+\frac{\partial H}{\partial\beta}\frac{\partial\beta}{\partial\varepsilon}\right)\times r\right] \\
&= E_{Y_0}\left[\frac{\partial H}{\partial\Delta}r+\frac{\sigma^2}{2}\frac{\partial^2 H}{\partial y_1^2}r\right] + E[\Delta_0^{-1}]E_{Y_0}\left[\frac{\partial H}{\partial\varepsilon}r\right]
\end{aligned}
$$

[recall that $\frac{\partial H}{\partial\beta}(Y_0,Y_0,0,\beta_0,0)=0$ under (19)]. This term may or may not be zero depending upon the functions $H$ and $r$. The next order term is given by

$$E_{\Delta,Y_0}\left[\Delta^{-1}\frac{\varepsilon^2}{2}(\Gamma^2_{\beta_0}\cdot(H\times r))(Y_0,Y_0,0,\beta_0,0)\right] = \frac{\varepsilon}{2}E_{\Delta,Y_0}[\Delta_0^{-1}(\Gamma^2_{\beta_0}\cdot(H\times r))].$$

Thus, plugging the result of (74) into (73), we get

$$T_\beta^H = \frac{2}{E[\Delta_0]}\bigg(\varepsilon^{\alpha-1}E_{\Delta,Y_0}[\Delta_0^{-1}(\Gamma_{\beta_0}\cdot H)\times r]$$
$$+ \frac{\varepsilon^\alpha}{2}E_{\Delta,Y_0}[\Delta_0^{-1}(\Gamma_{\beta_0}^2\cdot(H\times r))]\bigg) + O(\varepsilon^{\alpha+1}).$$

For multidimensional $H = (H_1,\ldots,H_r)'$, the $(i,j)$ term of the $T_\beta^H$ matrix is

$$[T_\beta^H]_{(i,j)} = \varepsilon^{\alpha_j-1}\frac{1}{E[\Delta_0]}E_{\Delta,Y_1,Y_0}[\Delta^{-1}H_i(Y_1,Y_0,\Delta,\bar\beta,\varepsilon)r_j(Y_1,\beta_0,\varepsilon)]$$
$$+ \varepsilon^{\alpha_i-1}\frac{1}{E[\Delta_0]}E_{\Delta,Y_1,Y_0}[\Delta^{-1}H_j(Y_1,Y_0,\Delta,\bar\beta,\varepsilon)r_i(Y_1,\beta_0,\varepsilon)]$$
$$+ O(\varepsilon^{\min(\alpha_i,\alpha_j)+1})$$
$$= \varepsilon^{\alpha_j-1}\frac{1}{E[\Delta_0]}E_{\Delta,Y_0}[\Delta_0^{-1}(\Gamma_{\beta_0}\cdot H_i)\times r_j]$$
$$+ \frac{\varepsilon^{\alpha_j}}{2}\frac{1}{E[\Delta_0]}E_{\Delta,Y_0}[\Delta_0^{-1}(\Gamma_{\beta_0}^2\cdot(H_i\times r_j))]$$
$$+ \varepsilon^{\alpha_i-1}\frac{1}{E[\Delta_0]}E_{\Delta,Y_0}[\Delta_0^{-1}(\Gamma_{\beta_0}\cdot H_j)\times r_i]$$
$$+ \frac{\varepsilon^{\alpha_i}}{2}\frac{1}{E[\Delta_0]}E_{\Delta,Y_0}[\Delta_0^{-1}(\Gamma_{\beta_0}^2\cdot(H_j\times r_i))]$$
$$+ O(\varepsilon^{\min(\alpha_i,\alpha_j)+1}).$$

6.6. *Proof of Theorem* 1. This corollary is a direct consequence of the (usual, nonstochastic) Taylor formula applied to the expression (14), with $D_\beta$ and $S_\beta$ given by Lemmas 1 and 3.

**7. Conclusions.** We have developed a set of tools for analyzing a large class of estimators of discretely-sampled continuous-time diffusions, including their asymptotic variance and bias. By Taylor-expanding the different matrices involved in the asymptotic distribution of the estimators, we are able to deliver fully explicit expressions of the various quantities determining the asymptotic properties of these estimators, and compare their relative merits. Our analysis covers the case where the sampling interval is random. As special cases, we cover the situation where the sampling is done at deterministic time-varying dates and the situation where the sampling occurs at fixed intervals. Most estimation methods can be analyzed within our framework—essentially any method that can be reduced to a method of moments or estimating equation problem. The two specific examples we analyzed display the various behaviors covered by our theorems, and we showed how our results can be used to assess the impact of different sampling patterns on the properties of these estimators.

Department of Economics
Princeton University and NBER
Princeton, New Jersey 08544-1021
USA
e-mail: yacine@princeton.edu

Department of Statistics
University of Chicago
Chicago, Illinois 60637-1514
USA
e-mail: mykland@galton.uchicago.edu